\documentclass[12pt]{amsart} 
\usepackage{amssymb,amscd}
\usepackage{verbatim}

\begin{document}

\numberwithin{equation}{section}

\newtheorem{theorem}[equation]{Theorem}
\newtheorem{lemma}[equation]{Lemma}
\newtheorem{conjecture}[equation]{Conjecture}
\newtheorem{proposition}[equation]{Proposition}
\newtheorem{corollary}[equation]{Corollary}
\newtheorem{cor}[equation]{Corollary}
\newtheorem*{theorem*}{Theorem}

\theoremstyle{definition}
\newtheorem*{definition}{Definition}
\newtheorem{example}[equation]{Example}

\theoremstyle{remark}
\newtheorem{remark}[equation]{Remark}
\newtheorem{remarks}[equation]{Remarks}
\newtheorem*{acknowledgement}{Acknowledgments}


\newenvironment{notation}[0]{%
  \begin{list}%
    {}%
    {\setlength{\itemindent}{0pt}
     \setlength{\labelwidth}{4\parindent}
     \setlength{\labelsep}{\parindent}
     \setlength{\leftmargin}{5\parindent}
     \setlength{\itemsep}{0pt}
     }%
   }%
  {\end{list}}

\newenvironment{parts}[0]{%
  \begin{list}{}%
    {\setlength{\itemindent}{0pt}
     \setlength{\labelwidth}{1.5\parindent}
     \setlength{\labelsep}{.5\parindent}
     \setlength{\leftmargin}{2\parindent}
     \setlength{\itemsep}{0pt}
     }%
   }%
  {\end{list}}
\newcommand{\Part}[1]{\item[\upshape#1]}

\renewcommand{\a}{\alpha}
\renewcommand{\b}{\beta}
\newcommand{\g}{\gamma}
\renewcommand{\d}{\delta}
\newcommand{\e}{\epsilon}
\newcommand{\f}{\phi}
\renewcommand{\l}{\lambda}
\renewcommand{\k}{\kappa}
\newcommand{\lhat}{\hat\lambda}
\newcommand{\m}{\mu}
\renewcommand{\o}{\omega}
\renewcommand{\r}{\rho}
\newcommand{\rbar}{{\bar\rho}}
\newcommand{\s}{\sigma}
\newcommand{\sbar}{{\bar\sigma}}
\renewcommand{\t}{\tau}
\newcommand{\z}{\zeta}

\newcommand{\D}{\Delta}

\newcommand{\gp}{{\mathfrak{p}}}
\newcommand{\gP}{{\mathfrak{P}}}
\newcommand{\gq}{{\mathfrak{q}}}
\newcommand{\gf}{{\mathfrak{f}}}

\newcommand{\Acal}{{\mathcal A}}
\newcommand{\Bcal}{{\mathcal B}}
\newcommand{\Ccal}{{\mathcal C}}
\newcommand{\Dcal}{{\mathcal D}}
\newcommand{\Ecal}{{\mathcal E}}
\newcommand{\F}{{\mathcal F}}
\newcommand{\Fcal}{{\mathcal F}}
\newcommand{\cF}{{\mathcal F}}
\newcommand{\Gcal}{{\mathcal G}}
\newcommand{\Hcal}{{\mathcal H}}
\newcommand{\Ical}{{\mathcal I}}
\newcommand{\Jcal}{{\mathcal J}}
\newcommand{\Kcal}{{\mathcal K}}
\newcommand{\Lcal}{{\mathcal L}}
\newcommand{\Mcal}{{\mathcal M}}
\newcommand{\Ncal}{{\mathcal N}}
\newcommand{\Ocal}{{\mathcal O}}
\newcommand{\Pcal}{{\mathcal P}}
\newcommand{\Qcal}{{\mathcal Q}}
\newcommand{\Rcal}{{\mathcal R}}
\newcommand{\Scal}{{\mathcal S}}
\newcommand{\Tcal}{{\mathcal T}}
\newcommand{\Ucal}{{\mathcal U}}
\newcommand{\Vcal}{{\mathcal V}}
\newcommand{\Wcal}{{\mathcal W}}
\newcommand{\Xcal}{{\mathcal X}}
\newcommand{\Ycal}{{\mathcal Y}}
\newcommand{\Zcal}{{\mathcal Z}}

\newcommand{\OO}{{\mathcal O}}    
\newcommand{\KK}{{\mathcal K}}    
\renewcommand{\O}{{\mathcal O}}   

\renewcommand{\AA}{\mathbb{A}}
\newcommand{\B}{\mathbb{B}}
\newcommand{\BB}{\mathbb{B}}
\newcommand{\CC}{\mathbb{C}}
\newcommand{\FF}{\mathbb{F}}
\newcommand{\GG}{\mathbb{G}}
\newcommand{\PP}{\mathbb{P}}
\newcommand{\NN}{\mathbb{N}}
\newcommand{\QQ}{\mathbb{Q}}
\newcommand{\RR}{\mathbb{R}}
\newcommand{\ZZ}{\mathbb{Z}}

\newcommand{\bfa}{{\mathbf a}}
\newcommand{\bfb}{{\mathbf b}}
\newcommand{\bfc}{{\mathbf c}}
\newcommand{\bfe}{{\mathbf e}}
\newcommand{\bff}{{\mathbf f}}
\newcommand{\bfg}{{\mathbf g}}
\newcommand{\bfp}{{\mathbf p}}
\newcommand{\bfr}{{\mathbf r}}
\newcommand{\bfs}{{\mathbf s}}
\newcommand{\bft}{{\mathbf t}}
\newcommand{\bfu}{{\mathbf u}}
\newcommand{\bfv}{{\mathbf v}}
\newcommand{\bfw}{{\mathbf w}}
\newcommand{\bfx}{{\mathbf x}}
\newcommand{\bfy}{{\mathbf y}}
\newcommand{\bfz}{{\mathbf z}}
\newcommand{\bfA}{{\mathbf A}}
\newcommand{\bfB}{{\mathbf B}}
\newcommand{\bfC}{{\mathbf C}}
\newcommand{\bfF}{{\mathbf F}}
\newcommand{\bfG}{{\mathbf G}}
\newcommand{\bfI}{{\mathbf I}}
\newcommand{\bfM}{{\mathbf M}}
\newcommand{\bfzero}{{\boldsymbol{0}}}
\newcommand{\bfmu}{{\boldsymbol\mu}}

\def\ta{{\tilde{a}}}
\def\tA{{\tilde{A}}}
\def\tb{{\tilde{b}}}
\def\tc{{\tilde{c}}}
\def\td{{\tilde{d}}}
\def\tf{{\tilde{f}}}
\def\tF{{\tilde{F}}}
\def\tg{{\tilde{g}}}
\def\tG{{\tilde{G}}}
\def\th{{\tilde{h}}}
\def\tK{{\tilde{K}}}
\def\tk{{\tilde{k}}}
\def\tz{{\tilde{z}}}
\def\tw{{\tilde{w}}}
\def\tphi{{\tilde{\varphi}}}
\def\tdelta{{\tilde{\delta}}}
\def\tbeta{{\tilde{\beta}}}
\def\talpha{{\tilde{\alpha}}}
\def\ttheta{{\tilde{\theta}}}
\def\tB{{\tilde{B}}}
\def\tR{{\tilde{R}}}
\def\tT{{\tilde{T}}}
\def\tE{{\tilde{E}}}
\def\tU{{\tilde{U}}}
\def\tGamma{{\tilde{\Gamma}}}

\newcommand{\ab}{{\textup{ab}}}
\newcommand{\Ahat}{{\hat A}}
\newcommand{\Aut}{\operatorname{Aut}}
\newcommand{\Cond}{{\mathfrak{N}}} 
\newcommand{\Disc}{\operatorname{Disc}}
\newcommand{\Div}{\operatorname{Div}}
\newcommand{\End}{\operatorname{End}}
\newcommand{\Frob}{\operatorname{Frob}}
\newcommand{\Fpbar}{{\overline{\FF}_p}}
\newcommand{\GK}{G_{\Kbar/K}}
\newcommand{\GL}{\operatorname{GL}}
\newcommand{\Gal}{\operatorname{Gal}}
\newcommand{\hhat}{{\hat h}}
\newcommand{\Image}{\operatorname{Image}}
\newcommand{\into}{\hookrightarrow}     
\newcommand{\Kbar}{{\bar K}}
\newcommand{\Lbar}{{\bar L}}
\newcommand{\ellbar}{{\bar \ell}}
\newcommand{\Kvbar}{{{\bar K}_v}}
\newcommand{\MOD}[1]{~(\textup{mod}~#1)}
\newcommand{\Norm}{\operatorname{N}}
\newcommand{\notdivide}{\nmid}
\newcommand{\nr}{{\textup{nr}}}    
\newcommand{\ord}{\operatorname{ord}}
\newcommand{\Pic}{\operatorname{Pic}}
\newcommand{\Qbar}{{\bar{\QQ}}}
\newcommand{\Qpbar}{{\bar{\QQ}_p}}
\newcommand{\kbar}{{\bar{k}}}
\newcommand{\QQbar}{{\bar{\QQ}}}
\newcommand{\rank}{\operatorname{rank}}
\newcommand{\res}{\operatornamewithlimits{res}}
\newcommand{\Resultant}{\operatorname{Res}}
\newcommand{\Res}{\operatorname{Res}}
\renewcommand{\setminus}{\smallsetminus}
\newcommand{\Spec}{\operatorname{Spec}}
\newcommand{\PGL}{\operatorname{PGL}}
\newcommand{\supp}{\operatorname{Supp}}
\newcommand{\tors}{{\textup{tors}}}
\newcommand{\val}{{\operatorname{val}}}
\newcommand{\<}{\langle}
\newcommand{\la}{{\langle}}
\renewcommand{\>}{\rangle}
\newcommand{\ra}{{\rangle}}
\newcommand{\Berk}{{\rm Berk}}
\newcommand{\BDV}{{\rm BDV}}
\newcommand{\Gauss}{{0}}
\newcommand{\Rat}{{\rm Rat}}
\newcommand{\RL}{{\rm RL}}
\newcommand{\HH}{{\mathbf H}}
\newcommand{\Gm}{{\mathbf G}_m}
\newcommand{\Hhat}{{\hat{H}}}


\hyphenation{archi-me-dean}


\title[Canonical heights over function fields]{A finiteness theorem for canonical heights attached to rational maps over function fields}

\author[Matthew Baker]{Matthew Baker}
\email{mbaker@math.gatech.edu}
\address{School of Mathematics,
          Georgia Institute of Technology, Atlanta GA 30332-0160, USA}

\date{First draft December 2005, revised May 2006}

\begin{abstract}
Let $K$ be a function field, 
let $\varphi \in K(T)$ be a rational map of degree $d \geq 2$ defined over $K$, 
and suppose that $\varphi$ is not isotrivial.
In this paper, we show that a point $P \in \PP^1(\Kbar)$ has $\varphi$-canonical height zero if and only
if $P$ is preperiodic for $\varphi$.  This answers affirmatively a question of Szpiro and Tucker, and generalizes a recent result
of Benedetto from polynomials to rational functions.  We actually prove the following stronger result, which is a variant of 
the Northcott finiteness principle: there exists $\varepsilon > 0$ such that the set of points $P \in \PP^1(K)$ with
$\varphi$-canonical height at most $\varepsilon$ is finite.
Our proof is essentially analytic, making use of potential theory on Berkovich spaces
to study the dynamical Green's functions $g_{\varphi,v}(x,y)$ attached to 
$\varphi$ at each place $v$ of $K$.  
For example, we show that every conjugate of $\varphi$ has 
bad reduction at $v$ if and only if $g_{\varphi,v}(x,x) > 0$ for all $x \in \PP^1_{\Berk,v}$, where
$\PP^1_{\Berk,v}$ denotes the Berkovich projective line over the completion of $\Kvbar$.
In an appendix, we use a similar method to give a new proof of the
Mordell-Weil theorem for elliptic curves over $K$.
\end{abstract}


\thanks{The author's research was supported by NSF Research Grant DMS-0300784.
The author would like to thank Rob Benedetto for helpful discussions related to this work,
Xander Faber for his many useful comments for improving the exposition, 
and Joe Silverman for his feedback on an earlier draft.}

\maketitle


\section{Introduction}
\label{IntroSection}

\subsection{Terminology}

Throughout this paper, $K$ will denote a {\em function field}, by which we will mean
a field endowed with a set $M_K$ of non-trivial non-archimedean absolute values which
satisfies the {\em product formula}:
\begin{itemize}
\item[(PF)] For every nonzero $x \in K$, $|x|_v = 1$ for all but finitely many $v \in M_K$ and
\[
\prod_{v \in M_K} |x|_v = 1 \ .
\]
\end{itemize}

Examples of function fields include the field of rational functions on any normal 
(or just regular in codimension 1) projective variety over a field $k$ 
(see \cite{LangDG},\S2.3, and \cite{BombieriGubler},\S1.4.6).



The {\em field of constants} of a function field $K$ is defined to be the subfield $k_0 \subset K$ consisting of all $x \in K$ such that
$|x|_v \leq 1$ for all $v \in M_K$.  By the product formula, if $x \in k_0$ is nonzero then in fact $|x|_v = 1$
for all $v \in M_K$.

\subsection{The canonical height attached to a rational map}

Let $\varphi \in K(T)$ be a rational map of degree $d \geq 2$ defined over $K$, so that
$\varphi$ acts on $\PP^1(K)$ in the usual way.
We define a {\em homogeneous lifting} of $\varphi$ to be a choice of
homogeneous polynomials $F_1,F_2 \in K[X,Y]$ of degree $d \geq 2$
having no common linear factor in $\Kbar[X,Y]$ such that
\[
\varphi([z_0:z_1]) \ = \ [F_1(z_0,z_1):F_2(z_0,z_1)]
\]
for all $[z_0:z_1] \in \PP^1(K)$.
(The polynomials $F_1,F_2$ are uniquely determined by $\varphi$ up to
multiplication by a common scalar $c \in K^*$.)
The mapping
\[
F = (F_1,F_2) : K^2 \to K^2
\]
is a lifting of $\varphi$ to $K^2$, and we denote by $F^{(n)} : K^2 \to K^2$ the
iterated map $F \circ F \circ \ldots \circ F$ ($n$ times).

The condition that $F_1$ and $F_2$ have no common linear factor over $\Kbar$ can be
rephrased by saying that the homogeneous resultant $\Res(F) = \Res(F_1,F_2)$ is nonzero.
If $F_1(z) = \prod_{i} z \wedge \alpha_i$ and 
$F_2(z) = \prod_{j} z \wedge \beta_j$, then
$\Res(F_1,F_2) = \prod_{i,j} \alpha_i \wedge \beta_j$, where
$(x_0,x_1) \wedge (y_0,y_1) = x_0 y_1 - x_1 y_0$.

\medskip

The canonical height function $\hat{H}_{F,K}$
is defined for $z \in K^2 \setminus \{ 0 \}$ by
\begin{equation*}
\begin{aligned}
\Hhat_{F,K}(P) &= \lim_{n\to\infty} \frac{1}{d^n} \sum_{v \in M_K} \log \| F^{(n)}(P) \|_v \ , \\
\end{aligned}
\end{equation*}
where $\| (x,y) \|_v = \max \{ |x|_v, |y|_v \}$.
It is straightforward to show using the product formula that 
$\hat{H}_{F,K}(z) \geq 0$, that $\hat{H}_{F,K}(z)$ depends only on the class of $z$ in
$\PP^1(K)$, and that $\hat{H}_{cF,K} = \hat{H}_{F,K}$ for all $c \in K^*$.
$\Hhat_{F,K}$ therefore descends
to a well-defined global canonical height function
$\hhat_{\varphi,K} : \PP^1(K) \to \RR_{\geq 0}$ depending only on $\varphi$.
For all $P \in \PP^1(K)$, $\hhat_{\varphi,K}$ satisfies the functional equation
\begin{equation}
\label{FunctEq}
\hhat_{\varphi,K}(\varphi(P)) = d \hhat_{\varphi,K}(P) \ .
\end{equation}

\medskip

If $K'/K$ is a finite extension, then for each place $w$ of $K'$ extending a given place $v$ of $K$, one
can define an absolute value on $K'$ via $|x|_w = |N^{K'_w}_{K_v}(x)|_v$.
The resulting set $M_{K'}$ of absolute values on $K'$ again satisfies (PF), and
for $P \in \PP^1(K)$ we have
\[
\hhat_{\varphi,K'}(P) = [K' : K] \hhat_{\varphi, K}(P) \ .
\]
(See \cite{SerreLMW},\S2.2, and \cite{BombieriGubler}, Proposition 1.4.2.)
It follows that we can extend $\hhat_{\varphi,K}$ to $\Kbar$ in a natural way by setting
$\hhat_{\varphi,K}(P) = \frac{1}{[K' : K]} \hhat_{\varphi, K'}(P)$
for $P \in \PP^1(K')$.  We will write $\hhat_{\varphi}(P)$ instead of $\hhat_{\varphi,K}(P)$
when no confusion is possible about which ground field we are talking about.

\medskip

If $h_K : \PP^1(K) \to \RR_{\geq 0}$ denotes the standard Weil height on $\PP^1(\Kbar)$, defined
for $P = (z_0:z_1) \in \PP^1(K)$ by
\[
h_K(P) = \sum_{v \in M_{K}} \log \max \{ |z_0|_v, |z_1|_v \}
\]
and extended to $\Kbar$ as above,
then there is a constant $C>0$ such that 
\begin{equation}
\label{BoundedDifference}
|\hhat_{\varphi,K}(P) - h_K(P)| \leq C \textrm{ for all } P \in \PP^1(\Kbar) \ .
\end{equation}
This follows from (\ref{FunctEq}) and the easily verified fact that
\begin{equation*}
\hhat_{\varphi,K}(P) = \lim_{n \to \infty} \frac{1}{d^n} h_K(\varphi^{(n)}(P))
\end{equation*}
for all $P \in \PP^1(\Kbar)$.

\medskip

It is easy to see using (\ref{FunctEq}) that $\hhat_{\varphi,K}(P) = 0$ if $P$ is preperiodic for $\varphi$
(i.e., if the orbit of $P$ under iteration of $\varphi$ is finite).
And if the field of constants $k_0$ of $K$ is finite,
then just as in the number field case, it is easy to show that $\hhat_{\varphi}$ satisfies
the following {\em Northcott finiteness property}: For any $M>0$, the set
\begin{equation}
\label{Northcott}
\{ P \in \PP^1(K) \; : \; \hhat_{\varphi,K}(P) \leq M \}
\end{equation}
is finite.
If property (\ref{Northcott}) holds for $K$, one easily deduces that if $\hhat_{\varphi,K}(P)=0$,
then $P$ is preperiodic.

However, in the function field case, we have:

\begin{lemma}
\label{NorthcottFailsLemma}
Let $K$ be a function field.  If the field of constants $k_0$ of $K$ is infinite, then
for $M$ sufficiently large, the set
\begin{equation}
\label{NorthcottFails}
\{ P \in \PP^1(K) \; : \; \hhat_{\varphi,K}(P) \leq M \}
\end{equation}
is infinite.
\end{lemma}

\begin{proof}
If $P \in \PP^1(k_0)$ then $h_K(P) = 0$, and therefore
$\hhat_{\varphi,K}(P) \leq C$ for all $P \in \PP^1(k_0)$ by (\ref{BoundedDifference}).
In particular, the set $\{ P \in \PP^1(K) \; : \; \hhat_{\varphi,K}(P) \leq C \}$ is infinite.
\end{proof}

\medskip

In addition, if $\varphi$ is defined over $k_0$, 
then $\hhat_{\varphi,K}(P) = 0$ for all $P \in \PP^1(k_0)$, but not
all points of $\PP^1(k_0)$ are preperiodic.  For example,
if $k_0 = \Qbar, K = \Qbar(T)$, and $\varphi(T) = T^2$,
then $\hhat_{\varphi,K}(P) = 0$ for all $P \in \PP^1(\Qbar)$, but the only preperiodic points in
$\PP^1(K)$ (or $\PP^1(\Kbar)$) are $0,\infty$, and the roots of unity in $\Qbar$.

More generally, if $\varphi$ is conjugate over a finite extension $K'$ of $K$ to a map defined over 
the field of constants of $K'$, then (\ref{Northcott}) fails over $K'$, since if 
$\varphi' = M^{-1} \circ \varphi \circ M$ with $M \in \PGL_2(K')$, then
\[
\hhat_{\varphi',K'}(P) = \hhat_{\varphi,K'}(M(P)) 
\]
for all $P \in \PP^1(K')$ (see Lemma~\ref{MobiusLemma} below).

\medskip

We will show, however, that a weak version of (\ref{Northcott}) still holds over function fields, 
even when the field of constants is infinite.  To state the result, 
define $\varphi \in K(T)$ to be {\em isotrivial over K} if 
there is a M{\"o}bius transformation $M \in \PGL_2(K)$ such that
$\varphi' = M^{-1} \circ \varphi \circ M$ is defined over the field of constants of $K$,
and {\em isotrivial} if there exists a finite extension 
$K'$ of $K$ such that $\varphi$ is isotrivial over $K'$.

\begin{theorem}
\label{maintheorem}
Let $K$ be a function field, and let $\varphi \in K(T)$ be a rational map of degree $d \geq 2$.
Assume that $\varphi$ is not isotrivial.  Then there exists $\varepsilon > 0$
(depending on $K$ and $\varphi$) such that the set
\[
\{ P \in \PP^1(K) \; : \; \hhat_{\varphi,K}(P) \leq \varepsilon \}
\]
is finite.
\end{theorem}

\begin{remark}
$\left.\right.$

(i) We do not know if Theorem~\ref{maintheorem} remains true if the hypothesis ``$\varphi$ is not
isotrivial'' is replaced by the {\em a priori} weaker hypothesis 
``$\varphi$ is not isotrivial {\em over K}.''

(ii) One cannot strengthen Theorem~\ref{maintheorem} to a Northcott-type theorem asserting,
under the same hypotheses as Theorem~\ref{maintheorem}, that
for each $D \geq 1$, there exists $\varepsilon > 0$ (depending on $D,K$, and $\varphi$) such that the set
\[
\{ P \in \PP^1(K') \; : \; [K' : K] \leq D \textrm{ and } \hhat_{\varphi,K}(P) \leq \varepsilon \}
\]
is finite.  For if this were true, then for all $n \geq 1$ the set
\[
\{ P \in \PP^1(K) \; : \; \hhat_{\varphi,K}(P) \leq \varepsilon d^n \}
\]
would be finite, since for each $P \in \PP^1(K)$ there exists $P'$ with $[K(P'):K] \leq d^n$ such that
$\varphi^{(n)}(P') = P$ and $\hhat_{\varphi,K}(P) = d^n \hhat_{\varphi,K}(P')$.
For $n$ sufficiently large, this contradicts Lemma~\ref{NorthcottFailsLemma}.

(iii) Combining Lemma~\ref{NorthcottFailsLemma} with Theorem~\ref{maintheorem}, we see that
if $k_0$ is infinite and $\varphi$ is not isotrivial, then the quantity 
$\varepsilon_{\varphi,K} := 
\inf_{P \in \PP^1(K)} \hhat_{\varphi,K}(P)$ is a positive real number.  It would be
interesting to investigate quantitatively the dependence of $\varepsilon_{\varphi,K}$ on
$\varphi$ and $K$. 
\end{remark}

\medskip

As a consequence of (\ref{FunctEq}) and Theorem~\ref{maintheorem}, 
we obtain:

\begin{cor}
\label{maincor}
If $K$ is a function field and $\varphi \in K(T)$ is a rational map of degree $d \geq 2$
which is not isotrivial, then a point $P \in \PP^1(\Kbar)$ satisfies $\hhat_{\varphi,K}(P)=0$
if and only if $P$ is preperiodic for $\varphi$.
\end{cor}

\begin{proof}
As previously mentioned, preperiodic points always have height zero.  Conversely, suppose 
$\hhat_{\varphi,K}=0$ and choose a finite extension $K'$ of $K$ with $P \in \PP^1(K')$.
Using (\ref{FunctEq}) and the fact that $\hhat_{\varphi,K'} = [K' : K] \hhat_{\varphi,K}$
for all $n \geq 0$, we obtain
\[
\hhat_{\varphi,K'}(\varphi^{(n)}(P))
= [K' : K] \hhat_{\varphi,K}(\varphi^{(n)}(P))
= d^n [K' : K] \hhat_{\varphi,K}(P) = 0 \ .
\]
By Theorem~\ref{maintheorem}, it follows that
the set $\{ P, \varphi(P), \varphi^{(2)}(P), \ldots \}$ is finite, i.e., 
that $P$ is preperiodic.
\end{proof}

The special case of Corollary~\ref{maincor} in which $\varphi$ is a {\em polynomial map} was proved
recently by Benedetto (\cite{BenedettoFFHeights}, Theorem B).  
In fact, Benedetto proves a more precise result in the polynomial case 
(\cite{BenedettoFFHeights}, Theorem A), as he just assumes that 
$\varphi$ is not isotrivial {\em over K}.
Corollary~\ref{maincor} answers affirmatively a question which we first
learned of from Lucien Szpiro and Thomas Tucker. 

\medskip

Our proof of Theorem~\ref{maintheorem}, which will be given in \S\ref{MainThmProofSection}, uses
the fact that isotriviality of $\varphi$
can be detected in terms of good and bad reduction.  In order to make this precise, we first recall the necessary
definitions.

Let $L$ be a valued field with valuation ring $\O_L$.  
If $\varphi \in L(T)$, we say that $\varphi$ has {\em good reduction} over $L$ if there 
exists a homogeneous lifting $F = (F_1,F_2)$ of $\varphi$ having good reduction, i.e.,
$F_1,F_2 \in \O_L[x,y]$ and $\Res(F_1,F_2) \in \O_L^*$.  
If $\varphi$ does not have good reduction over $L$, then we say that $\varphi$ has {\em bad reduction} over $L$.

Also, we say that $\varphi$ has {\em potentially good reduction} over $L$ if
there is a finite extension 
$L'$ of $L$ and a M{\"o}bius transformation $M \in \PGL_2(L')$ such that
$\varphi' = M^{-1} \circ \varphi \circ M$ has good reduction over $L'$.

Conversely, we say that $\varphi$ has {\em genuinely bad reduction} over $L$ if
it does not have potentially good reduction over $L$.

In \S\ref{IsotrivialitySection}, we will prove the following criterion for isotriviality:

\begin{theorem}
\label{IsotrivialityTheorem}
Let $K$ be a function field, and let $\varphi \in K(T)$ be a rational map of degree at least 2.
Then $\varphi$ is isotrivial if and only if $\varphi$ has potentially good reduction over $K_v$
for all $v \in M_K$, where $K_v$ denotes the completion of $K$ at $v$.
\end{theorem}

\subsection{Dynamical Green's functions and reduction}
\label{GFSection}

In order to apply Theorem~\ref{IsotrivialityTheorem} to canonical heights, we will use a decomposition
\begin{equation}
\label{eq:globaleq1}
\hhat_{\varphi,K}(x) + \hhat_{\varphi,K}(y) = \sum_{v \in M_K} g_{\varphi,v}(x,y)
\end{equation}
of the global canonical height into local contributions.
Here $g_{\varphi,v}(x,y)$ is a two-variable {\em dynamical Green's function} attached to 
$\varphi$ which was introduced in \cite{BREQUI}; see
\S\ref{DGFSection} for a definition and further details.
The function $g_{\varphi,v}(x,y)$ is defined 
first for $x$ and $y$ in
the completion $\CC_v$ of $\overline{K}_v$, and then
extended in a natural way to the Berkovich projective line
$\PP^1_{\Berk}$ over $\CC_v$ (see \S\ref{BerkovichSection}).

In \S\ref{ReductionTheoremSection}, we will use the theory of Berkovich spaces 
to prove a result characterizing
genuinely bad reduction in terms of dynamical Green's functions.  
In order to state the result, let $L$ be a complete and algebraically closed non-archimedean field, 
let $\varphi \in L(T)$ be a rational map of degree at least 2,
and let $g_{\varphi}(x,y)$ be the corresponding dynamical Green's function on $\PP^1(L)$.
We view $\PP^1(L)$ as being endowed with the analytic topology coming from the norm on $L$.

\begin{theorem}
\label{PositivityTheorem}
If $\varphi$ has genuinely bad reduction, then there exists a constant $\beta > 0$ and a covering of $\PP^1(L)$ by
finitely many analytic open sets $V_1,\ldots,V_s$ such that for each $1 \leq i \leq s$, we have
$g_\varphi(x,y) \geq \beta$ for all $x,y \in V_i$.
\end{theorem}

As we will see, Theorem~\ref{PositivityTheorem} can be easily deduced from 
the following result concerning the extension of $g_{\varphi}(x,y)$ to the Berkovich
projective line:

\begin{theorem}
\label{PositivityTheorem2}
$\varphi$ has genuinely bad reduction over $L$ if and only if 
$g_{\varphi}(x,x) > 0$ for all $x \in \PP^1_{\Berk} \backslash \PP^1(L)$.
\end{theorem}

We will derive Theorem~\ref{PositivityTheorem2} as a special case of a more general result 
(Theorem~\ref{ReductionTheorem}) giving 
several different conditions which are all equivalent to $\varphi$ having potentially good reduction.
It is worth noting that Theorem~\ref{PositivityTheorem2}, and therefore the theory of
Berkovich spaces, is a crucial ingredient in our proof of Theorem~\ref{maintheorem}, 
even though the statement of the latter result has nothing to do with Berkovich spaces.

\begin{remark}
The results in \cite{BRBook} are currently
stated and proved in the special case where $L=\CC_p$.  However, all of the results from \cite{BRBook} 
relevant to the present paper remain valid over an
arbitrary complete and algebraically closed non-archimedean field, see
\cite{BRBook2} 
and \cite{Thuillier} for details.
\end{remark}

\medskip

Another (this time non-essential) ingredient in our proof of Theorem~\ref{maintheorem} is the following 
lower bound for average values of $g_\varphi$ which was proved in \cite{BakerDGF}, Theorem 1.1:

\begin{theorem}
\label{GreenTheorem}
Let $L$ be a valued field, and let $\varphi \in L(T)$ be a rational map of degree $d \geq 2$.
Then there is an effective constant $\alpha> 0$, depending on $\varphi$ and $K$, such that if 
$N \geq 2$ and $z_1,\ldots,z_N$ are distinct points of $\PP^1(K)$, then
\begin{equation*}
\sum_{\substack{1 \leq i,j \leq N \\ i \neq j}} g_{\varphi}(z_i,z_j)
\geq -\alpha N \log N \ .
\end{equation*}
\end{theorem}

\medskip

\begin{remark}
$\left.\right.$

(i) As noted above, Corollary~\ref{maincor} is a generalization of Benedetto's main theorem in
\cite{BenedettoFFHeights} from polynomials to rational functions.  As with the passage from the main result
of \cite{BakerHsia} to that of \cite{BREQUI}, the rational function case 
of Corollary~\ref{maincor} seems to require more machinery than the polynomial case.  
For example, 
the analogue of Theorem~\ref{PositivityTheorem} was proved in \cite{BenedettoFFHeights} by 
analyzing the radii of preimages of suitable disks under $\varphi$.  However, we have not found a way 
to avoid using the machinery of Berkovich spaces to prove Theorem~\ref{PositivityTheorem} for rational maps.
Although the proofs of Theorems~\ref{IsotrivialityTheorem} and \ref{GreenTheorem} do not rely on Berkovich's theory,
they are also more technically
involved than their polynomial counterparts in \cite{BenedettoFFHeights}.

(ii) The interested reader should compare our results with those of 
Moriwaki \cite{Moriwaki}, who obtains a 
Northcott-type theorem for varieties over function fields by replacing the usual Weil height 
(which he calls a ``geometric height'') with an ``arithmetic height'' coming from
the choice of a certain kind of polarization.  

(iii) It would be interesting to formulate and prove a result analogous to Theorem~\ref{maintheorem} in
higher dimensions.
\end{remark}


\section{Proof of Theorem~\ref{maintheorem}}
\label{MainThmProofSection}

We now give the proof of our main theorem on canonical heights over function fields
(Theorem~\ref{maintheorem}),
assuming Theorems~\ref{IsotrivialityTheorem} and \ref{PositivityTheorem}.

\begin{proof}
(Compare with Theorem~1.14 of \cite{BakerDGF}.)

Recall from (\ref{eq:globaleq1})
that for all $z,w \in \PP^1(K)$ with $z\neq w$, we have
\[
\sum_{v \in M_{K}} g_{\varphi,v}(z,w) = \hhat_{\varphi}(z) + \hhat_{\varphi}(w) \ .
\]
Since $\varphi$ is not isotrivial, it follows from 
Theorem~\ref{IsotrivialityTheorem} that there exists a place $v_0 \in M_K$
at which $\varphi$ has genuinely bad reduction.
Let $\CC_{v_0}$ be the smallest complete and algebraically closed field containing $K_{v_0}$.
By Theorem~\ref{PositivityTheorem},
there is a finite covering $V_1,\ldots,V_s$ of $\PP^1(\CC_{v_0})$ and a constant $C_1 > 0$ 
such that $g_{\varphi,v_0}(z,w) \geq C_1$ whenever $z,w \in V_i$ ($i=1,\ldots,s$).
If $z_1,\ldots,z_N \in \PP^1(K)$, then 
by Theorem~\ref{GreenTheorem} and the fact that $g_{\varphi,v} \geq 0$ whenever 
$\varphi$ has good reduction at $v$,
there is a constant $C_2>0$ depending on $\varphi$ and $K$ such that
\[
g(z_1,\ldots,z_N) := \sum_{v \in M_K} \sum_{i \neq j} g_{\varphi,v}(z_i,z_j) \geq -C_2 N\log N \ .
\]
Moreover, if $z_1,\ldots,z_N \in V_i$ for some $i$, then
\begin{equation}
\label{eq:globaleq2}
g(z_1\ldots,z_N) \geq C_1 N^2 - C_2 N\log N \ .
\end{equation}

Let $M = [\frac{N-1}{s}]+1$.  By the pigeonhole principle, in any subset of $\PP^1(\CC_{v_0})$ of cardinality $N$, there is an $M$-element subset contained in
some $V_i$.  Without loss of generality, order the $z_j$'s so that this subset is $\{ z_1,\ldots,z_{M} \}$.  
Applying (\ref{eq:globaleq1}) and (\ref{eq:globaleq2}), we obtain
\[
C_1 M^2 - C_2 M\log M \leq g(z_1,\ldots,z_{M}) \leq 2M^2 \max_j \hhat_{\varphi}(z_j) \ .
\]

If $\hhat_{\varphi}(z_j) \leq \frac{C_1}{4}$ for all $j=1,\ldots,N$, then we obtain
\[
\frac{C_1}{2} M \leq C_2 \log M \ ,
\]
which implies that $M \leq C_3$ for some constant $C_3 > 0$ depending 
only on $\varphi$ and $K$.
Thus $N \leq Ms+1 \leq C_4$ for some $C_4>0$ depending only on $\varphi$ and $K$.
Setting $\varepsilon = \frac{C_1}{4}$ now gives the desired result.  
\end{proof}

\begin{remark}
In the proof of Theorem~\ref{maintheorem}, one could avoid appealing to Theorem~\ref{GreenTheorem}
by applying Theorem~\ref{PositivityTheorem} and a pigeonhole principle argument simultaneously at each place
of genuinely bad reduction for $\varphi$.  This would provide a better value of $\varepsilon$ but worse 
upper bound on $N$ than is given by the argument above.
\end{remark}


\section{Dynamical Green's functions, the Berkovich projective line, and reduction of rational maps}
\label{DGFSection}

In this section, we study the dynamical Green's function $g_{\varphi}(z,w)$ attached to a rational map $\varphi$ of degree $d \geq 2$
defined over a valued field $L$.  A proof of Theorem~\ref{PositivityTheorem} is given in \S\ref{ReductionTheoremSection}.

\subsection{Definition and basic properties of dynamical Green's functions}

We begin by recalling some terminology from \cite{BakerDGF} and \cite{BREQUI}.
As in \S\ref{IntroSection}, write $\varphi$ in the form
\[
\varphi([z_0:z_1]) \ = \ [F_1(z_0,z_1):F_2(z_0,z_1)]
\]
for some homogeneous polynomials $F_1,F_2 \in L[x,y]$ of degree $d \geq 2$
with $\Res(F_1,F_2) \neq 0$, and let 
$F = (F_1,F_2) : L^2 \to L^2$.

Write $\|(z_0,z_1)\| = \max \{ |z_0|,|z_1| \} $, and
for 
$z \in L^2 \backslash \{ 0 \}$,
define the {\em homogeneous local dynamical height} $\Hhat_{F} : L^2
\backslash \{ 0 \} \to \RR$ by
\[
\Hhat_{F}(z) = \lim_{n\to\infty} \frac{1}{d^n} \log \| F^{(n)}(z) \| \ .
\]

By convention, we set $\Hhat_{F}(0,0) = -\infty$.
According to \cite{BREQUI}, Lemma 3.5, the limit $\lim_{n\to\infty} \frac{1}{d^n} \log \| F^{(n)}(z) \|$
exists for all $z \in L^2 \backslash \{ 0 \}$, and 
$\frac{1}{d^n} \log \| F^{(n)}(z) \|$ converges uniformly on
$L^2 \backslash \{ 0 \}$ to $\Hhat_{F}(z)$.  

\medskip

The function $\Hhat_F$ {\em scales logarithmically}, in the sense that for every $\alpha \in L^*$ and
$z \in L^2 \backslash \{ 0 \}$,
\begin{equation}
\label{eq:logscale}
\Hhat_F(\alpha z) = \Hhat_F(z) + \log |\alpha| \ .
\end{equation}

If the value group $\{ -\log |z| \; : \; z \in L^* \}$ is dense in $\RR$
(which will be the case, for example, if $L$ is algebraically closed and $| \cdot |$ is non-trivial), 
it follows from (\ref{eq:logscale}) that given $\varepsilon > 0$, for every
$z \in L^2 \backslash \{ 0 \}$ there exists $\alpha \in L^*$ such that
\begin{equation}
\label{eq:logscale2}
-\varepsilon < \Hhat_F(\alpha z) < 0 \ .
\end{equation}



\medskip

When $z, w \in L^2$ are linearly independent over $L$, define
\[
G_{F}(z,w) \ = \ -\log |z\wedge w| + \Hhat_{F}(z) + \Hhat_{F}(w)
+ \log R \ ,
\]
where $R = |\Res(F)|^{-\frac{1}{d(d-1)}}$.

According to \cite{BREQUI}, Lemma 3.21, for all 
$\alpha, \beta,\gamma \in L^*$, we have
\[
G_{\gamma F}(\alpha z, \beta w) =  G_{F}(z,w) \ . 
\]

In particular, $G_{F}$ descends to a well-defined
function $g_{\varphi}(z,w)$ on $\PP^1(L)$:
for $z, w \in \PP^1(L)$ and any lifts $\tz, \tw \in L^2 \backslash \{ 0 \}$,
\begin{equation} 
\label{gFvDef} 
g_{\varphi}(z,w) \ = \ -\log |\tz\wedge \tw| + \Hhat_{F}(\tz) 
                      + \Hhat_{F}(\tw) + \log R \ .
\end{equation}
If $z\neq w$ then the right-hand side of (\ref{gFvDef}) is finite; if
$z=w$ then we set $g_{\varphi}(z,z) = +\infty$.


\medskip

If $K$ is a field satisfying (PF) and $\varphi \in K(T)$ is a rational map of degree $d \geq 2$, 
then for each $v \in M_K$ we have (setting $L = \CC_v$) 
an associated dynamical Green's function
$g_{\varphi,v}(z,w)$.
By the product formula (applied twice), 
for all $z,w \in \PP^1(K)$ with $z\neq w$ we have the fundamental identity
\begin{equation*}
\sum_{v \in M_K} g_{\varphi,v}(z,w) = \hhat_{\varphi}(z) + \hhat_{\varphi}(w) \ .
\end{equation*}

\subsection{Extending the dynamical Green's function to Berkovich space}
\label{BerkovichSection}

There is a natural extension of $g_\varphi(z,w)$ to 
the Berkovich projective line; we briefly recall 
the relevant background material from \cite{BRBook},\cite{BRBook2}.

Let $L$ be a complete non-archimedean field.  
With its usual topology, the projective line $\PP^1(L)$ is totally disconnected,
and if the residue field of $L$ is infinite then $\PP^1(L)$ is not locally compact.
The Berkovich projective line $\PP^1_{\Berk} = \PP^1_{\Berk , L}$ over $L$ is a 
connected compact Hausdorff space which contains $\PP^1(L)$ as a dense subspace.
The construction of $\PP^1_{\Berk}$
is functorial, and in particular a rational map $\varphi : \PP^1(L) \to \PP^1(L)$
extends in a natural way to a continuous map from $\PP^1_{\Berk}$ to itself.

The space $\PP^1_{\Berk}$ has the structure of a metric $\RR$-tree (in the sense of \cite{FJBook}),
and admits a rich {\em potential theory}, including a theory of harmonic and subharmonic functions.  
These notions are defined in terms of a measure-valued Laplacian operator $\Delta$ 
on $\PP^1_{\Berk}$ whose domain is the space $\BDV(\PP^1_{\Berk})$ of 
functions of {\em bounded differential variation}
(see \cite{BRBook}, \S5.3).  
A similar theory has been developed independently by at least three different sets of
authors -- see \cite{FJBook}, \cite{BRBook}, and \cite{Thuillier}.  For simplicity, 
we will use \cite{BRBook} and \cite{BRBook2} as our
basic reference for potential theory on $\PP^1_{\Berk}$, although the results we need are also 
contained in \cite{Thuillier}.

We recall from \cite{BREQUI} that an Arakelov Green's function on
the Berkovich projective line is a function
$g(z,w) : \PP^1_{\Berk} \times \PP^1_{\Berk}
\rightarrow \RR \cup \{ +\infty\}$
such that

\begin{itemize}
\item[(B1)] (Semicontinuity)
The function $g(z,w)$ is lower-semicontinuous, 
and is finite and continuous off the diagonal.

\item[(B2)] (Differential equation)  For each $w \in \PP^1_{\Berk}$,
$g(z,w)$ belongs to the space $\BDV(\PP^1_{\Berk})$.  Furthermore,
there is a probability measure $\mu$ on $\PP^1_{\Berk}$
such that for each $w$, $g(z,w)$ satisfies the identity
\[
\Delta_z g(z,w) \ = \ \delta_w(z) - \mu(z) \ .
\]
\end{itemize}

Conditions (B1) and (B2) imply that $g(z,w)$ is symmetric and bounded below
(see \cite{BRBook}, Proposition 7.19).
The semicontinuity along the diagonal
is a technical condition which arises naturally from properties of the space
$\PP^1_{\Berk}$ (see \cite{BRBook}, Proposition 3.1).
Together, (B1) and (B2) determine
$g(z,w)$ up to an additive constant by the maximum principle
(\cite{BRBook}, Proposition 5.14).  If in addition
\begin{itemize}
\item[(B3)] (Normalization) \quad $\iint_{\PP^1_{\Berk} \times \PP^1_{\Berk}} g(z,w) \, \mu(z) \mu(w) \ = \ 0$ \ ,
\end{itemize}
we will say $g(z,w)$
is a {\em normalized} Arakelov Green's function.

Since it is known that
\begin{equation}
\label{eq:GreenPotentialConstant}
\int_{\PP^1_{\Berk}} g(z,w) \mu(z) \equiv C
\end{equation}
for some constant $C$ independent of $w$ (see \cite{BRBook}, Proposition 7.24), 
(B3) is in fact equivalent to the {\em a priori} stronger condition
\begin{itemize}
\item[$({\rm B3})^\prime$]\quad $\int g(z,w) \, \mu(z) \ \equiv \ 0$ \ .
\end{itemize}

\medskip

If $\varphi \in L(T)$ is a rational function of degree $d \geq 2$, then
as explained in \cite{BREQUI}, \S7.5, the dynamical Green's function 
$g_\varphi$ extends in a natural way to a normalized Arakelov Green's function 
\begin{equation}
g_{{\varphi}}(z,w) : \PP^1_{\Berk} \times \PP^1_{\Berk} \to \RR \cup \{ +\infty \} \ .
\end{equation}
Explicitly, we have
\[
g_{\varphi}(z,w) = \liminf_{\substack{(z_0,w_0) \in \PP^1(L) \times \PP^1(L) \\ 
(z_0,w_0) \to (z,w)}} g_{\varphi}(z_0,w_0) \ .
\]
It is precisely because of the constant $\log R$ which appears in the definition of $g_\varphi$
that it satisfies (B3); see Appendix~\ref{Appendix1} for details.

We will refer to the probability measure $\mu_{\varphi}$ on $\PP^1_{\Berk}$ 
defined by the differential equation 
$\Delta_z g_{\varphi}(z,w) = \delta_w - \mu_{\varphi}$
as the {\em canonical measure} associated to $\varphi$.  
The fact that $\mu_{\varphi}$ is a {\em probability} measure (i.e., that $\mu_{\varphi}$ is 
nonnegative and has total mass 1) is proved in \cite{BRBook}, Theorem 7.14.
It is shown in \cite{BRBook}, Proposition 7.15 that the $\mu_{\varphi}$ 
has no point masses on $\PP^1(L)$.
However, it can have point masses on $\PP^1_{\Berk} \backslash
\PP^1(L)$ (see Theorem~\ref{ReductionTheorem}).

\begin{lemma}
\label{ContinuityLemma}
For each fixed $w \in \PP^1_{\Berk}$, the 
function $f_w(z) = g_\varphi(z,w)$ is continuous on all of $\PP^1_{\Berk}$ as an 
extended real-valued function.
If $w \in \PP^1_{\Berk} \backslash \PP^1(L)$, then $f_w(z)$ is real-valued on all of $\PP^1_{\Berk}$, and if $w \in \PP^1(L)$, then
$f_w(z)$ is real-valued except at $z=w$, where $f_w(w)=+\infty$.
\end{lemma}

\begin{proof}
See Proposition 7.19 of \cite{BRBook}.
\end{proof}

There is a notion of subharmonic functions on $\PP^1_{\Berk}$ (see \cite{BRBook}, \S6, or \cite{Thuillier}, \S3.1.2, for
a definition).  We will use the fact that subharmonic functions satisfy the following {\em maximum principle}:

\begin{theorem}[Maximum Principle for Subharmonic Functions]
\label{MaximumPrinciple}
If $U$ is an open subset of $\PP^1_{\Berk}$ and $f : U \to \RR \cup \{ +\infty \}$ is a subharmonic function which attains its
maximum value on $U$, then $f$ is locally constant.
\end{theorem}

\begin{proof}
See \cite{BRBook}, Proposition 6.15, and also
\cite{Thuillier}, Proposition 3.1.11.
\end{proof}

As a consequence, we can identify where the minimum and maximum values
of $g_{\varphi}$ occur:

\begin{theorem}
\label{MinMaxTheorem}
For each fixed $w \in \PP^1_{\Berk}$, the 
minimum value of $f_w(z) = g_\varphi(z,w)$ 
is achieved on $\supp(\mu_{\varphi})$, and the maximum value is 
achieved at $w$.
\end{theorem}

\begin{proof}
By Lemma~\ref{ContinuityLemma}, $f_w(z)$ attains its maximum and minimum values on
$\PP^1_{\Berk}$.  (The maximum value is $+\infty$ if $w \in \PP^1(L)$.)
By Proposition 6.1, Definition 6.2, and Proposition 7.19 of \cite{BRBook}, 
$f_w(z)$ is subharmonic on $\PP^1_{\Berk} \backslash \{ w \}$ and
$-f_w(z)$ is subharmonic on the complement of $\supp(\mu_{\varphi})$.
(This can also be deduced from \cite{Thuillier}, Proposition 3.4.4.)
The result therefore follows from the maximum principle for subharmonic functions 
(Theorem~\ref{MaximumPrinciple}).
\end{proof}

\begin{remark}
Theorem~\ref{MinMaxTheorem} 
holds more generally for Arakelov-Green's functions associated to ``log-continuous'' 
probability measures, see \cite{BRBook}, \S7.
\end{remark}

\subsection{Homogeneous filled Julia sets and transfinite diameters}
\label{HomogFJSection}

Let $L$ be a complete valued field, and as in the previous section, let $F = (F_1,F_2)$ be a non-degenerate 
homogeneous polynomial mapping of $L^2$, i.e.,
$F_1,F_2 \in L[x,y]$ are homogeneous polynomials of degree $d \geq 2$
such that $\Res(F) \neq 0$.  

The {\emph{homogeneous filled Julia set}} $K_{F}$ of $F$ in $L^2$
is the set of all $z \in L^2$ for which $\| F^{(n)}(z) \|$
stays bounded as $n$ tends to infinity.
Clearly $F(K_F) = K_F$ and $F^{-1}(K_{F}) = K_{F}$, i.e., $K_F$ is completely invariant under $F$.
According to \cite{BREQUI}, Lemma 3.8, 
\[
K_{F} = \{ z \in L^2 \; : \; \Hhat_{F}(z) \leq 0 \} \ .
\]


If $F$ has good reduction, then by
\cite{BREQUI}, Lemma 3.9,
$K_{F} = B(0,1)$ is the unit polydisc in
$L^2$ and $\Hhat_{F}(z) = \log \| z \|$ for all $z \in L^2$.
We will establish a converse to this result shortly (Corollary~\ref{KFBallCor}).

Also, when $F$ has good reduction, the dynamical Green's function is given quite 
simply by
\begin{equation}
\label{GoodRedGreen}
G_F(z,w) = -\log|z \wedge w| + \log \| z \| + \log \| w \| \ .
\end{equation}
Note that the right-hand side of (\ref{GoodRedGreen}) is always non-negative.
This means that $g_\varphi(x,y) \geq 0$ for all $x,y \in \PP^1(L)$ when $\varphi$ has good
reduction.  By Lemma~\ref{MobiusLemma}, the same is true when $\varphi$ has potentially
good reduction.
Conversely, we will see in Theorem~\ref{ReductionTheorem} that if
$g_\varphi(x,y) \geq 0$ for all $x,y \in \PP^1(L)$, then $\varphi$ has potentially good
reduction.

\medskip

Recall that for $z = (z_1,z_2), w = (w_1,w_2) \in L^2$, we have
$z \wedge w = z_1 w_2 - z_2 w_1$.  
By analogy with the classical transfinite diameter,
if $E \subset L^2$ is a bounded set, we define
\[
d^0_n(E) \ = \ \sup_{z_1,\ldots,z_n \in E}
\left( \prod_{i \neq j} |z_i \wedge z_j|  \right)^{\frac{1}{n(n-1)}} 
\]

By \cite{BREQUI}, Lemma 3.10, the sequence of nonnegative real numbers $d^0_n(E)$ is
non-increasing.  In particular, the quantity $d^0_\infty(E) =
\lim_{n\to\infty} d^0_n(E)$ is well-defined.
We call $d^0_\infty(E)$ the {\em homogeneous transfinite diameter} of
$E$.

We also define $d^0(E)$ to be 
\[
d^0(E) = \sup_{z,w \in E} |z \wedge w| \ .
\]
We call $d^0(E)$ the {\em wedge diameter} of $E$.

As an example, it is easy to see that if $L$ is a complete and algebraically closed non-archimedean field
and $r \in |L^*|$, then the wedge diameter and homogeneous transfinite diameter of the polydisc 
\[
B(0,r) = \{ z \in L^2 \; : \; \| z \| \leq r \}
\]
are both equal to $r^2$.

If $M \in \GL_2(L)$, then it is elementary to show that $|M(z) \wedge M(w)| = |\det(M)| |z \wedge w|$,
and thus $d^0_\infty(M(E)) = |\det(M)| d^0_\infty(E)$ and 
$d^0(M(E)) = |\det(M)| d^0(E)$.

\medskip

As explained in \cite{BakerDGF} (see also \cite{BREQUI}), for any valued field $L$ we have
\begin{equation}
\label{TransfiniteDiameterInequality}
d^0_\infty(K_{F}) \ \leq \ |\Res(F)|^{-1/d(d-1)} \ .
\end{equation}
If $L$ is complete and algebraically closed, then 
equality holds in (\ref{TransfiniteDiameterInequality})
by Corollary~\ref{TransfiniteDiameterCor}.



\begin{lemma}
\label{IntegralityLemma}
If $L$ is a complete and algebraically closed non-archimedean field and $F(B(0,1)) \subseteq B(0,1)$, 
then all coefficients of $F_1$ and $F_2$ lie in $\O_L$.
\end{lemma}

\begin{proof}
Write $F(x,y) = (a_0 x^d + a_1 x^{d-1}y + \cdots + a_d y^d,
b_0 x^d + b_1 x^{d-1}y + \cdots + b_d y^d)$.  
Since $F(B(0,1)) \subseteq B(0,1)$, in particular 
the one-variable polynomials $a(x) = a_0 x^d + a_1 x^{d-1} + \cdots + a_d$ and 
$b(x) = b_0 x^d + b_1 x^{d-1} + \cdots + b_d$ map $\O_L$ to itself.  
But then we are done, since
\[
\sup_{z \in \O_L} |a(z)| = \max_{0 \leq i \leq d} |a_i| \ ,
\]
and similarly for $\sup_{z \in \O_L} |b(z)|$.
\end{proof}

\begin{cor}
\label{KFBallCor}
If $L$ is a complete and algebraically closed non-archimedean field and $K_F = B(0,1)$,
then $F$ has good reduction.
\end{cor}

\begin{proof}
Since $K_F = B(0,1)$ implies that $F(B(0,1)) \subseteq B(0,1)$, it follows from 
Lemma~\ref{IntegralityLemma} that all coefficients of $F$ lie in $\O_L$.
Also, $1 = d^0_\infty(K_F) = |\Res(F)|^{-\frac{1}{d(d-1)}}$, so $|\Res(F)|=1$.
It follows that $F$ has good reduction.
\end{proof}

\subsection{Dynamical Green's functions and reduction of rational maps}
\label{ReductionTheoremSection}

The main result of this section is the following theorem:

\begin{theorem}
\label{ReductionTheorem}
Let $L$ be a complete and algebraically closed non-archimedean field, and let $\varphi \in L(T)$ be a 
rational map of degree $d \geq 2$.  Then $g_{\varphi}(x,x) \geq 0$ for all $x \in \PP^1_{\Berk}$, and 
the following are equivalent:
\begin{itemize}
\item[(1)] $\varphi$ has potentially good reduction.
\item[(2a)] $g_{\varphi}(x,y) \geq 0$ for all $x,y \in \PP^1_{\Berk}$.
\item[(2b)] $g_{\varphi}(x,y) \geq 0$ for all $x,y \in \PP^1(L)$.
\item[(3a)] $g_{\varphi}(\zeta,\zeta) = 0$ for some $\zeta \in \PP^1_{\Berk}$.
\item[(3b)] $g_{\varphi}(\zeta,\zeta) = 0$ for some $\zeta \in \PP^1_{\Berk} \backslash \PP^1(L)$.
\item[(4a)] $\mu_{\varphi} = \delta_{\zeta}$ for some $\zeta \in \PP^1_{\Berk}$.
\item[(4b)] $\mu_{\varphi} = \delta_{\zeta}$ for some $\zeta \in \PP^1_{\Berk} \backslash \PP^1(L)$.
\end{itemize}
\end{theorem}

Note that conditions (2a) and (2b) are equivalent by the lower semicontinuity of 
$g_{\varphi}(x,y)$, (3a) and (3b) are equivalent because $g_\varphi(x,x)=+\infty$ for 
all $x \in \PP^1(L)$, and (4a) and (4b) are equivalent because, as mentioned in \S2.2, 
$\mu_{\varphi}$ has no point masses in $\PP^1(L)$.

As an immediate corollary of Theorem~\ref{ReductionTheorem}, we have:

\begin{cor}
\label{ReductionCor}
Under the same hypotheses as Theorem~\ref{ReductionTheorem}, the following are equivalent:
\begin{itemize}
\item[(1)] $\varphi$ has genuinely bad reduction.
\item[(2)] $g_{\varphi}(x,y) < 0$ for some $x,y \in \PP^1(L)$.
\item[(3)] $g_{\varphi}(\zeta,\zeta) > 0$ for all $\zeta \in \PP^1_{\Berk} \backslash \PP^1(L)$.
\item[(4)] $\mu_{\varphi}$ is not a point mass.
\end{itemize}
\end{cor}

Assuming Corollary~\ref{ReductionCor}, we show how to deduce 
Theorem~\ref{PositivityTheorem}, which is needed for the proof of Theorem~\ref{maintheorem}.  

\begin{proof}[Proof of Theorem~\ref{PositivityTheorem}]
(Compare with \cite{BakerDGF}, Lemma~3.1.)

Suppose $\varphi$ has genuinely bad reduction over $L$.
We claim that there exists a constant $\beta > 0$ such that 
$g_{\varphi}(x,x) > \beta$ for all $x \in \PP^1_{\Berk}$.
Indeed, according to Corollary~\ref{ReductionCor}, 
the fact that $\varphi$ has genuinely bad reduction means that
$g_{\varphi}(x,x) > 0$ for all $x \in \PP^1_{\Berk}$.
The claim then follows from the fact that $\PP^1_{\Berk}$ is compact and 
$g_{\varphi}$ is lower semicontinuous (see \cite{BourbakiGT}, \S~IV.6.2, Theorem 3).
The lower semicontinuity of $g_{\varphi}(z,w)$, together with 
the definition of the product topology on $\PP^1_{\Berk} \times \PP^1_{\Berk}$,
also shows that for each $x \in \PP^1_{\Berk}$, there is an open neighborhood
$U_x$ of $x$ in $\PP^1_{\Berk}$ such that $g_{\varphi}(z,w) > \beta$ for $z,w \in U_x$.
By compactness of $\PP^1_{\Berk}$, the covering $\{ U_x \; | \; x \in \PP^1_{\Berk} \}$ has a finite
subcover $U_1,\ldots,U_s$.  Take $V_i = U_i \cap \PP^1(L)$.

\end{proof}

Before giving the proof of Theorem~\ref{ReductionTheorem}, we establish some useful lemmas.

\begin{lemma}
\label{GoodReductionLiftLemma}
If $L$ is an algebraically closed non-archimedean field and $\varphi$
has potentially good reduction over $L$, then for any homogeneous lifting $F$ of $\varphi$, 
there exists $M \in \GL_2(L)$ such that $M^{-1} \circ F \circ M$ has good reduction.
\end{lemma}

\begin{proof}
Choose any $M' \in \GL_2(L)$ such that $\varphi' = M'^{-1} \circ \varphi \circ M'$ has good reduction.
This means that there exists a homogeneous lifting $F'$ of $\varphi'$ with good reduction.  Since any two lifts 
of $\varphi'$ differ by a nonzero constant, we have $F' = c(M'^{-1} \circ F \circ M')$ for some $c \in L^*$.  It is enough to
prove that there exists $\alpha \in L^*$ such that $F' = M^{-1} \circ F \circ M$, where $M = \alpha M'$.  
For this, we take $\alpha$ such that $\alpha^{d-1} = c$, and then compute that
\[
\begin{aligned}
(M^{-1} \circ F \circ M)(z) &= \alpha^{-1} (M'^{-1} \circ F \circ M')(\alpha z) \\
&= \alpha^{d-1}(M'^{-1} \circ F \circ M')(z) \\
&= F'(z) \\
\end{aligned}
\]
as desired.
\end{proof}

\begin{lemma}
\label{LiftingPreperiodicPointsLemma}
Let $L$ be an algebraically closed field, let $\varphi \in L(T)$ be a rational map of degree $d \geq 2$, and 
let $F = (F_1,F_2)$ be a homogeneous lifting of $\varphi$.
Then the natural map from the set of preperiodic points of $F$ in $L^2 \backslash \{ 0 \}$ to the set of preperiodic points of $\varphi$
in $\PP^1(L)$ is surjective.
\end{lemma}

\begin{proof}
Let $P \in \PP^1(L)$ be a preperiodic point of $\varphi$, so that $\varphi^{(j)}(P) = \varphi^{(k)}(P)$ for
some positive integers $j \neq k$.  Then for any homogeneous lifting $z \in L^2 \backslash \{ 0 \}$ of $P$,
there exists a nonzero constant $c \in L$ such that $F^{(j)}(z) = c F^{(k)}(z)$.
Choose $\alpha \in L^*$ such that $\alpha^{d^k} = c \alpha^{d^j}$.  
Then 
\[
F^{(j)}(\alpha z) = \alpha^{d^j}F^{(j)}(z) = c \alpha^{d^j} F^{(k)}(z)
= \alpha^{d^k} F^{(k)}(z) = F^{(k)}(\alpha z) \ ,
\]
so that $\alpha z \in L^2$ is a lifting of $P$ which is preperiodic for $F$.
\end{proof}

\begin{lemma}
\label{NormalFormLemma}
Let $L$ be an algebraically closed field, and let $\varphi \in L(T)$ be a rational map of degree $d \geq 2$.
Let $w,w' \in \PP^1(L)$ be distinct preperiodic points of $\varphi$.
Then there exists a M{\"o}bius transformation $M \in \PGL_2(L)$ with
$M(\infty)=w$ and $M(0)=w'$, and
a lifting $F' = (F_1',F_2')$ of $\varphi' = M^{-1} \circ \varphi \circ M$ such that
$(1,0)$ and $(0,1)$ are preperiodic points for $F'$.
(Here we identify the point $(1:0) \in \PP^1(L)$ with $\infty$ and
$(0:1) \in \PP^1(L)$ with $0$.)
\end{lemma}

\begin{proof}
Choose any $M'' \in \GL_2(L)$ such that $M''(\infty)=w$ and $M''(0)=w'$, so
that $0,\infty$ are preperiodic for $\varphi'' = M''^{-1} \circ \varphi \circ M''$.  
Choose any homogeneous lifting $F'' = (F_1,F_2)$ of $\varphi''$ to $L^2$.
By Lemma~\ref{LiftingPreperiodicPointsLemma}, there are
lifts $(0,s),(t,0) \in L^2$ of $0$ and $\infty$, respectively, such that
$(0,s)$ and $(t,0)$ are preperiodic for $F''$.  
Let $F' = M'^{-1} \circ F'' \circ M'$, where $M' \in \GL_2(L)$ is the linear transformation
such that $M'(x,y) = (tx,sy)$.  Then $(0,1)$ and $(1,0)$ are preperiodic for $F'$,
the M{\"o}bius transformation $M = M'' \circ M'$ satisfies
$M(\infty)=w$ and $M(0)=w'$, and $F'$ is a lifting of $\varphi' = M^{-1} \circ \varphi \circ M$ 
as desired.
\end{proof}

\begin{lemma}
\label{PreperiodicGreenLemma}
Let $L$ be an algebraically closed field, and let $\varphi \in L(T)$ be a rational map of degree $d \geq 2$.
If $g_{\varphi}(x,y) \geq 0$ for all $x,y \in \PP^1(L)$, then there exist preperiodic points 
$w,w' \in \PP^1(L)$ such that $g_{\varphi}(w,w')=0$.
\end{lemma}

\begin{proof}
Let $z \in \PP^1(L)$ be a fixed point of $\varphi$, and let $w \in \PP^1(L)$ 
be any periodic point of $\varphi$ different from
$z$.  Let $N$ be the period of $w$.  By Theorem~\ref{PullbackTheorem}, if 
\[
(\varphi^{(N)})^*(z) = z + \sum_{i=1}^{d^N - 1} z_i' 
\]
as divisors on $\PP^1(L)$, then 
\begin{equation}
\label{eq:PeriodicGreenSum}
g_{\varphi}(w,z) = g_{\varphi}(\varphi^{(N)}(w),z)
= g_{\varphi}(w,z) + \sum_{i=1}^{d^N - 1} g_{\varphi}(w,z_i') \ .
\end{equation}

Since $w \neq z$, all terms appearing in (\ref{eq:PeriodicGreenSum}) are finite, and
therefore 
\[
\sum_{i=1}^{d^N - 1} g_{\varphi}(w,z_i') = 0 \ .
\]
Since $g_{\varphi}(x,y) \geq 0$ for all $x,y \in \PP^1(L)$ by assumption, it follows that
$g_{\varphi}(w,z_i') = 0$ for all $i$.  
Since $\varphi^{(N)}(z_i')=z$, each $z_i'$ is preperiodic, so we can 
take $w' = z_1'$.  
\end{proof}

\medskip

We are now ready to prove Theorem~\ref{ReductionTheorem}.

\begin{proof}[Proof of Theorem~\ref{ReductionTheorem}]

For the first assertion, 
we know that $g_{\varphi}(x,x)=+\infty$ for $x \in \PP^1(L)$, so it
suffices to show that if $\zeta \in \PP^1_{\Berk} \backslash \PP^1(L)$, then
$g_{\varphi}(\zeta,\zeta) \geq 0$.
By Theorem~\ref{MinMaxTheorem}, for all $x \in \PP^1_{\Berk}$ we have
\[
g_{\varphi}(\zeta,\zeta) \geq g_{\varphi}(x,\zeta) \ .
\]
Integrating against the probability measure $\mu$ with respect to the variable $x$, we obtain
\[
g_{\varphi}(\zeta,\zeta) \geq \int_{\PP^1_{\Berk}} g_{\varphi}(x,\zeta) \mu(x) = 0 
\]
by $({\rm B3})^\prime$.


Now for the various equivalences.  
We have already remarked that $(2a) \Leftrightarrow (2b),(3a) \Leftrightarrow (3b),(4a) \Leftrightarrow (4b)$.
Let us therefore refer to these assertions as $(2),(3),(4)$, respectively.
We will show that $(3) \Leftrightarrow (4)$, $(4) \Leftrightarrow (2)$, and finally that 
$(1) \Rightarrow (2)$ and $(4) \Rightarrow (1)$.

\medskip

$(3) \Rightarrow (4)$: 

Suppose $g_{\varphi}(\zeta,\zeta) = 0$ for some $\zeta \in \PP^1_{\Berk} \backslash \PP^1(L)$.
Then by Theorem~\ref{MinMaxTheorem}, for all $x \in \PP^1_{\Berk}$, we have 
\begin{equation}
\label{eq:gneg}
0 = g_{\varphi}(\zeta,\zeta) \geq g_{\varphi}(x,\zeta) \ .
\end{equation}
On the other hand, since $\int_{\PP^1_{\Berk}} g_{\varphi}(x,\zeta) \mu(x) = 0$, it
follows from Lemma~\ref{ContinuityLemma} that $g_{\varphi}(x,\zeta) = 0$ for all 
$x \in \supp(\mu_{\varphi})$.  
Since (by Theorem~\ref{MinMaxTheorem}) the minimum value of $g_{\varphi}(x,\zeta)$ is achieved
on $\supp(\mu_{\varphi})$, it follows from ({eq:gneg}) that 
$g_{\varphi}(x,\zeta) = 0$ for all $x \in \PP^1_{\Berk}$.
But then 
\[
\delta_{\zeta} - \mu_{\varphi} = \Delta_x g_{\varphi}(x,\zeta) = 0 
\]
and thus $\mu_{\varphi} = \delta_{\zeta}$.


\medskip

$(4) \Rightarrow (3)$: If $\mu_{\varphi} = \delta_{\zeta}$, then since $\mu_{\varphi}$ is
normalized (Corollary~\ref{NormalizedCor}), we have
\[
0 = \iint g_{{\varphi}}(x,y) \mu_{\varphi}(x) \mu_{\varphi}(y) = 
g_{{\varphi}}(\zeta,\zeta) \ .
\]

\medskip

$(4) \Rightarrow (2)$: If $\mu_{\varphi} = \delta_{\zeta}$, then
\[
g_{{\varphi}}(x,\zeta) = 
\int g_{{\varphi}}(x,y) \mu_{\varphi}(y) = 0
\]
for all $x \in \PP^1_{\Berk}$ by $({\rm B3})^\prime$.

Fix $y \in \PP^1_{\Berk}$, and let $f_y(x) = 
g_{{\varphi}}(x,y)$.
By Theorem~\ref{MinMaxTheorem},
the minimum value of $f_y(x)$ on $\PP^1_{\Berk}$ occurs at $\zeta$, and 
therefore for all $x \in \PP^1_{\Berk}$, 
\[
g_{{\varphi}}(x,y) \geq g_{{\varphi}}(x,\zeta) = 0 \ .
\]

\medskip

$(2) \Rightarrow (4)$: 

Since 
\[
\int g_{{\varphi}}(x,y) \mu_{\varphi}(y) = 0
\]
for all $x \in \PP^1_{\Berk}$,
if $g_{\varphi}(x,y) \geq 0$ for all $x,y \in \PP^1_{\Berk}$
and $\zeta \in \supp(\mu_{\varphi})$, we must have $g_{\varphi}(\cdot,\zeta) \equiv 0$
by Lemma~\ref{ContinuityLemma}.
Therefore
\[
\delta_{\zeta} - \mu_{\varphi} = \Delta_x g_{\varphi}(x,\zeta) = 0 \ ,
\]
which means that $\mu_{\varphi} = \delta_{\zeta}$.

\medskip

$(1) \Rightarrow (2)$: 

If $\varphi$ has good reduction, 
this follows from (\ref{GoodRedGreen}).  The 
case of potentially good reduction then follows from
Lemma~\ref{MobiusLemma}.

\medskip

$(4) \Rightarrow (1)$: 

Assuming that $\mu_{\varphi} = \delta_\zeta $,
we want to show that $\varphi$ has potentially good reduction.

Since $(4) \Rightarrow (2)$, 
Lemma~\ref{PreperiodicGreenLemma} implies that we can find preperiodic points
$w,w' \in \PP^1(L)$ for $\varphi$ such that $g_{\varphi}(w,w') = 0$.  Replacing 
$\varphi$ by a conjugate if necessary, Lemma~\ref{NormalFormLemma} shows that we can assume that
$w = 0, w' = \infty$, and that there is a homogeneous lifting $F$ of $\varphi$ so that
$(1,0)$ and $(0,1)$ are preperiodic for $F$.  In particular, 
we have $\Hhat_F((0,1)) = \Hhat_F((1,0)) = 0$.
We will show that $F$ has good reduction.  
By Corollary~\ref{KFBallCor}, it suffices to prove that $K_F = B(0,1)$.

Recall that for any lifts $z,z' \in L^2 \backslash \{ 0 \}$ of $P,P' \in \PP^1(L)$, 
respectively, we have
\[
g_\varphi(P,P') = G_F(z,z') = -\log |z \wedge z'| + \Hhat_F(z) + \Hhat_F(z') + \log |\Res(F)|^{-\frac{1}{d(d-1)}} \ .
\]
In particular,
\[
0 = g_\varphi(w,w') = G_F((0,1),(1,0)) = \log |\Res(F)|^{-\frac{1}{d(d-1)}} \ ,
\]
so that $|\Res(F)|=1$ and 
\[
G_F(z,z') = -\log |z \wedge z'| + \Hhat_F(z) + \Hhat_F(z')
\]
for all $z,z' \in L^2 \backslash \{ 0 \}$.

Take $z = (x,y) \in K_F$.  Then $\Hhat_F(z) \leq 0$, so 
condition (2) implies that
\[
0 \leq G_F(z,(0,1)) = -\log |x| + \Hhat_F(z) \leq -\log |x|
\]
and thus $|x| \leq 1$.  An identical calculation using $(1,0)$ 
shows that $|y| \leq 1$ as well.  Thus $z \in B(0,1)$, and we
conclude that $K_F \subseteq B(0,1)$.

Let $\varepsilon > 0$, and 
take $z = (x,y) \in B(0,1)$ and $z' = (x',y') \in K_F$
so that $-\varepsilon < \Hhat_F(z') < 0$.
Since $z,z' \in B(0,1)$, we have $|z \wedge z'| \leq 1$, and thus
\begin{equation}
\label{eq:GHinequality}
G_F(z,z') \geq \Hhat_F(z) - \varepsilon \ .
\end{equation}

Let $P$ be the image of $z$ in $\PP^1(L)$, and let $P' \in \PP^1(L)$ be arbitrary.
By (\ref{eq:logscale2}), we may choose a lift $z' \in L^2 \backslash \{ 0 \}$ of $P'$ such that 
$-\varepsilon < \Hhat_F(z') < 0$.
Since $g_{\varphi}(x,y)$ is continuous off the diagonal
and $\mu_{\varphi}= \delta_{\zeta}$ implies that 
$g_{\varphi}(\zeta,y) \equiv 0$, if $P' \to \zeta$ then 
\[
G_F(z',z) = g_{\varphi}(P',P) \to 
g_{\varphi}(\zeta,P) = 0.
\]
As $\PP^1(L)$ is dense in $\PP^1_{\Berk}$, we may use 
(\ref{eq:GHinequality}) to conclude that
$\Hhat_F(z) \leq \varepsilon$.
Letting $\varepsilon \to 0$ shows that $\Hhat_F(z) \leq 0$.
Thus $z \in K_F$ and $K_F = B(0,1)$ as desired.
\end{proof}

\begin{remark}
The fact that $g_{\varphi}(x,x) \geq 0$ for all $x \in \PP^1_{\Berk}$, as well as
the implication $(3) \Rightarrow (4)$ in the proof of 
Theorem~\ref{ReductionTheorem}, can also be deduced from the following
``Energy Minimization Principle'' (see \cite{BRBook}, Theorem 7.20):

Define the ``energy functional'' $I_\varphi(\nu)$ on the space $\PP$ of probability
measures on $\PP^1_{\Berk}$ by the formula
\[
I_\varphi(\nu) \ = \
\iint_{\PP^1_{\Berk} \times \PP^1_{\Berk}} g_\varphi(z,w) \, \nu(z) \nu(w) \ .
\]
Then $I_\varphi(\nu) \geq I_\varphi(\mu_{\varphi}) = 0$ for all $\nu \in \PP$,
with equality if and only if $\nu = \mu_{\varphi}$.
\end{remark}

As a consequence of Theorem~\ref{ReductionTheorem},
we present new proofs of some theorems of R.~Benedetto and 
J.~Rivera-Letelier concerning reduction of rational maps.
For example, we easily deduce from 
Theorem~\ref{ReductionTheorem} the following result of Benedetto 
(\cite{BenedettoRDJ}, Theorem B, see also \cite{RLPHS}, \S7):

\begin{cor}
\label{BenedettoCor}
$\varphi$ has potentially good reduction iff $\varphi^{(n)}$ has potentially good reduction for some (equivalently, every)
$n \geq 1$.
\end{cor}

\begin{proof}
This follows from Theorem~\ref{ReductionTheorem},
together with the fact that 
\begin{equation*}
g_{\varphi}(z,w) = g_{\varphi^{(n)}}(z,w)
\end{equation*}
for any $n \geq 1$,
which follows easily from 
the definition of $g_{\varphi}(z,w)$ using the fact that
\[
|\Res(F)|^{-\frac{1}{d(d-1)}} = 
|\Res(F^{(n)})|^{-\frac{1}{d^n(d^n-1)}}
\]
by (\ref{ResultantComp2}).
\end{proof}

The {\em exceptional set} $\Ecal(\varphi)$ of $\varphi$ in $\PP^1_{\Berk} \backslash \PP^1(L)$,
is defined to be the set of all $\zeta \in \PP^1_{\Berk} \backslash \PP^1(L)$ such that 
\[
\bigcup_{n \geq 0} \varphi^{(-n)}(\zeta)
\]
is finite.
The following result is originally due to Rivera-Letelier (see \cite{RLPHS}, \S7, Theorem 3):

\begin{cor}
If $\varphi$ has genuinely bad reduction, then $\Ecal(\varphi) = \emptyset$.
If $\varphi$ has potentially good reduction, then
$\Ecal(\varphi) = \{ \zeta \}$ consists of a single point.
\end{cor}

\begin{proof}
Since $\varphi$ maps $\PP^1_{\Berk} \backslash \PP^1(L)$ surjectively onto itself,
the forward orbit of any point in $\Ecal(\varphi)$ must be a finite cycle of length $n$.  
Using Corollary~\ref{BenedettoCor}, after replacing $\varphi$ by $\varphi^{(n)}$
it suffices to prove that there is at most one point $\zeta \in 
\PP^1_{\Berk} \backslash \PP^1(L)$ which is completely invariant under $\varphi$, 
and that there exists a completely invariant point if and only if $\varphi$ has potentially 
good reduction.  

If $\zeta \in \PP^1_{\Berk} \backslash \PP^1(L)$ is completely invariant, then
$\varphi^*(\zeta) = d\zeta$ as divisors on $\PP^1_{\Berk}$, 
and by Theorem~\ref{PullbackTheorem} we have
\[
g_\varphi(\zeta,\zeta) = g_\varphi(\varphi(\zeta),\zeta)
= g_\varphi(\zeta,\varphi^*(\zeta)) = d \cdot g_\varphi(\zeta,\zeta) \ ,
\]
so that $g_\varphi(\zeta,\zeta) = 0$.  By Theorem~\ref{ReductionTheorem}, this implies
that $\varphi$ has potentially good reduction, and that $\mu_{\varphi} = \delta_\zeta$.
In particular, there can be at most one such point $\zeta$. 
Conversely, if $\varphi$ has potentially good reduction, then 
$\mu_{\varphi} = \delta_\zeta$ for some $\zeta \in \PP^1_{\Berk} \backslash \PP^1(L)$.
Since $\varphi^*(\mu_{\varphi}) = d\mu_{\varphi}$ and $\varphi_*(\mu_{\varphi}) = \mu_{\varphi}$, 
it follows easily that $\zeta$ is completely invariant under $\varphi$.
\end{proof}

\begin{remark}
Let $\zeta_{\Gauss}$ denote the {\em Gauss point} of $\PP^1_{\Berk}$ 
(see \cite{BRBook}, \S1).
Then it is not hard using the above
considerations to show the equivalence of the following conditions:
\begin{itemize}
\item[(1)] $\varphi$ has good reduction.
\item[(2)] $\varphi^{(n)}$ has good reduction for some (equivalently, every)
$n \geq 1$.
\item[(3)] $\varphi^{-1}(\zeta_{\Gauss}) = \{ \zeta_\Gauss \}$.
\item[(4)] $\mu_{\varphi} = \delta_{\zeta_{\Gauss}}$.
\end{itemize}
We leave the details to the reader.
The equivalence of (1) and (2) is originally due to Benedetto
(\cite{BenedettoRDJ}, Theorem B), and 
the equivalence of (1) and (3) was first noted by Rivera-Letelier in 
\cite{RLPHS}, \S4.1.  Rivera-Letelier also shows in \cite{RLPHS} that
$\varphi$ has {\em non-constant} reduction if and only if
$\varphi(\zeta_{\Gauss}) = \zeta_\Gauss$.
\end{remark}


\section{A criterion for isotriviality}
\label{IsotrivialitySection}

We recall the statement of Theorem~\ref{IsotrivialityTheorem}:

\begin{theorem*}
Let $K$ be a function field, and let $\varphi \in K(T)$ be a rational map of degree at least 2.
Then $\varphi$ is isotrivial if and only if $\varphi$ has potentially good reduction over $K_v$
for all $v \in M_K$.
\end{theorem*}

\begin{proof}
It follows from the definitions that isotrivial maps have everywhere potentially good reduction.
We therefore assume that $\varphi$ has potentially good reduction over $K_v$ for all $v$, and want to
show that $\varphi$ is isotrivial.

Since $\PP^1(\Kbar)$ contains infinitely many preperiodic points for $\varphi$,
after conjugating $\varphi$ by a M{\"o}bius transformation and replacing $K$ by a finite extension if necessary,
we may assume without loss of generality by Lemma~\ref{NormalFormLemma} that
there is a homogeneous lifting $F$ of $\varphi$ for which $(0,1)$ and $(1,0)$ are preperiodic.
In particular, $P = (0,1)$ and $Q = (1,0)$ belong to the homogeneous filled Julia set $K_{F,v}$ for all
$v \in M_K$.

Let $v \in M_K$ be a (non-archimedean) place of $K$, let $\CC_v$ be the smallest complete and
algebraically closed field containing the completion $K_v$ of $v$, and let $\OO_v$ be its valuation ring.  
Let $E_v$ be the homogeneous filled Julia set of $F$ in $\CC_v^2$.
By Lemma~\ref{GoodReductionLiftLemma},
there exists $M = M_v \in \GL_2(\CC_v)$ such that
$F' = M^{-1} \circ F \circ M$ has good reduction over $\CC_v$.
Write $M^{-1}(P)=P'$ and $M^{-1}(Q)=Q'$; since $E'_v = M^{-1}(E_v)$ is the homogeneous
filled Julia set for $F'$, we have $E'_v = B(0,1)$ in $\CC_v^2$.
As $aP' + bQ' \in B(0,1)$ for all $a,b \in \CC_v$ with $|a|,|b| \leq 1$,
and since $M$ is linear and takes $B(0,1)$ to $E_v$, it follows that $aP + bQ \in E_v$ for 
all $|a|,|b| \leq 1$.  Since $P = (0,1)$ and $Q = (1,0)$, we thus have
$B(0,1) \subseteq E_v$.  

In particular, $d^0_\infty(E_v) \geq 1$ for all $v \in M_K$.
But according to (\ref{TransfiniteDiameterInequality}) and the product formula,
we have 
\[
\prod_{v \in M_K} d^0_\infty(E_v) \leq 1 \ .
\]
Therefore $d^0_\infty(E_v) = 1$ for all $v \in M_K$.

Since $E_v' = B(0,1)$, letting $d^0$ denote the ``wedge diameter'' of a set 
as in \S\ref{HomogFJSection}, we have
\[
d^0_\infty(E_v) =   |\det(M)| d^0_\infty(E_v')
= |\det(M)| d^0(E_v') = d^0(E_v)
\]
for all $v \in M_K$.  
Thus $B(0,1) \subseteq E_v$ and $d^0(E_v) = 1$
for all $v \in M_K$.  
But this implies that $E_v = B(0,1)$ for all $v \in M_K$, since 
if $(x,y) \in E_v$, then $|x| = |(x,y) \wedge (0,1)| \leq 1$ and
$|y| = |(x,y) \wedge (1,0)| \leq 1$, so that $(x,y) \in B(0,1)$.

To conclude, 
write $F(x,y) = (a_0 x^d + a_1 x^{d-1}y + \cdots + a_d y^d,
b_0 x^d + b_1 x^{d-1}y + \cdots + b_d y^d)$ and
note that $E_v = B(0,1)$ implies that $F(B(0,1)) \subseteq B(0,1)$, and hence 
$F_1,F_2 \in \O_v[x,y]$ by Lemma~\ref{IntegralityLemma}.
But then $|a_i|_v,|b_i|_v \leq 1$ for all $0 \leq i \leq d$ and all $v \in M_K$,
which means that all coefficients of $F$ lie in the constant field $k_0$ of $K$.
Thus $\varphi$ is isotrivial as desired.
\end{proof}

\appendix
\section{Functional equation for the dynamical Green's function and consequences}
\label{Appendix1}

Let $L$ be a valued field, 
and let $\varphi \in L(T)$ be a 
rational map of degree $d \geq 2$.  
Our goal in this appendix is to prove a functional equation for 
the dynamical Green's function $g_{\varphi}(z,w)$ (Theorem~\ref{PullbackTheorem}),
and to deduce some consequences of this formula.  Some of these results are also proved in \cite{BREQUI} 
in the case where $L$ is the completion of a global field.  The proofs
given here, by contrast, are purely local.

\medskip

Let $k$ be a field.  By 
\cite{LangAlgebra}, Theorem IX.3.13,
if $F = (F_1,F_2)$ and $G = (G_1,G_2)$ 
where $F_1,F_2 \in k[X,Y]$ are homogeneous of degree $d$ and
$G_1,G_2 \in k[X,Y]$ are homogeneous of degree $e$, then
\begin{equation}
\label{ResultantComp}
\Res(F \circ G) = \Res(F)^e \Res(G)^{d^2}.
\end{equation}

In particular, it follows by induction on $n$ that
\begin{equation}
\label{ResultantComp2}
\Res(F^{(n)}) = \Res(F)^{\frac{d^{n-1}(d^n - 1)}{d-1}} \ .
\end{equation}

\begin{lemma}
\label{MobiusLemma}
Let $L$ be a valued field, and let $\varphi \in L(T)$ be a 
rational map of degree $d \geq 2$.  Let $F$ be a homogeneous lifting of $\varphi$, let $M \in \GL_2(L)$,
and let $F' = M^{-1}\circ F \circ M$.
Then for all $z,w \in L^2 \backslash \{ 0 \}$,
\begin{equation}
\label{MHeight}
\Hhat_F(M(z)) = \Hhat_{F'}(z)
\end{equation}
and
\begin{equation}
\label{MGreen}
G_F(M(z),M(w)) = G_{F'}(z,w) \ .
\end{equation}
\end{lemma}

\begin{proof}
First of all note that given $M$, there exist constants $C_1,C_2 > 0$ such that
\begin{equation}
\label{MSandwich}
C_1 \| z \| \leq \| M(z) \| \leq C_2 \| z \| \ .
\end{equation}
Indeed, if we take $C_2$ to be the maximum of the absolute values of the entries of $M$ then
clearly $\| M(z) \| \leq C_2 \| z \|$ for all $z$.  
By the same reasoning, if we let $C_1^{-1}$ be
the maximum of the absolute values of the entries of $M^{-1}$, then
\[
\| M^{-1} (M (z)) \| \leq C_1^{-1} \| M(z) \| \ ,
\]
which gives the other inequality.

By the definition of $\Hhat_F$, we have 
\[
\begin{aligned}
\Hhat_F(M(z)) & = \lim_{n \to \infty} \frac{1}{d^n} \log \| F^{(n)}(M(z)) \| \\
& = \lim_{n \to \infty} \frac{1}{d^n} \log \| M F'^{(n)}(z) \|
= \Hhat_{F'}(z) \ , \\
\end{aligned}
\]
where the last equality follows from (\ref{MSandwich}).
This proves (\ref{MHeight}).

Since $|M(z) \wedge M(w)| = |\det(M)| \cdot |z \wedge w|$, we have
\begin{equation}
\label{eq:MWedge}
-\log | M(z) \wedge M(w)| = -\log |\det(M)| - \log |z \wedge w| \ .
\end{equation}

On the other hand, by (\ref{ResultantComp}) and the fact that $|\Res(M)| = |\det(M)|$, 
it follows that
\begin{equation}
\label{eq:ResTilde}
-\frac{1}{d(d-1)} \log |\Res(F')| = 
-\frac{1}{d(d-1)} \log |\Res(F)| - \log |\det(M)| \ .
\end{equation}

Putting together (\ref{MHeight}),(\ref{eq:MWedge}), and
(\ref{eq:ResTilde}) gives (\ref{MGreen}). 
\end{proof}

\begin{cor}
\label{MobiusCor}
Let $L$ be a valued field, and let $\varphi \in L(T)$ be a 
rational map of degree $d \geq 2$.  Let $M \in \PGL_2(L)$.  
Then for all $z,w \in \PP^1(L)$, we have
\[
g_\varphi(M(z),M(w)) = g_{M^{-1}\circ \varphi \circ M}(z,w) \ . 
\]
\end{cor}

The following result gives a useful functional equation for the dynamical Green's function.
A proof using the fact that $g_{\varphi}$ is {\em normalized} (i.e., satisfies condition (B3) above)
is given in \cite{BREQUI}, Corollary 3.39.  However, the proof in \cite{BREQUI} that $g_{\varphi}$ is
normalized uses global methods.  Here, we give a simple, purely local proof of the functional equation, and
instead deduce from the functional equation that $g_{\varphi}$ is normalized (Corollary~\ref{NormalizedCor}).

\begin{theorem}
\label{PullbackTheorem}
Let $L$ be an algebraically closed valued field, and let $\varphi \in L(T)$ be a 
rational map of degree $d \geq 2$.  
Then for all $x,y \in \PP^1(L)$, we have
\[
g_\varphi(\varphi(x),y) = g_\varphi(x,\varphi^{*}(y)) \ , 
\]
where $y_1,\ldots,y_d$ are the preimages of $y$ under $\varphi$, counting
multiplicity, and
\[
g_\varphi(x,\varphi^{*}(y)) := \sum_{i=1}^d
g_\varphi(x,y_i) \ . 
\]
\end{theorem}

\begin{proof}
Using Corollary~\ref{MobiusCor}, we may assume without loss of 
generality that $y = \infty$.
Let $F$ be a homogeneous lifting of $\varphi$, and let
$R = |\Res(F)|^{-\frac{1}{d(d-1)}}$.  For $z \in L^2 \backslash \{ 0 \}$, let 
$[z]$ denote the class of $z$ in $\PP^1(L)$.  Fix
$w = (1,0)$, let $a_1, \ldots a_d$ be the preimages of $[w]$ under $\varphi$ (counting
multiplicities), and let $w_1,\ldots,w_d$ be solutions to $F(w_i) = w$ such that $[w_i] = a_i$.

It suffices to show that for all $z \in L^2 \backslash \{ 0 \}$, we have
\begin{equation}
\label{eq:PullbackId1}
G_F(F(z),w) = \sum_{i=1}^d G_F(z,w_i) \ .
\end{equation}
Since $\Hhat_F(F(z)) = d \Hhat_F(z)$ and 
$\Hhat_F(w_i) = \frac{1}{d} \Hhat_F(w)$ for all $i$, (\ref{eq:PullbackId1}) is
equivalent to
\begin{equation*}
-\log |F(z) \wedge w| = (d-1)\log R - \sum_{i=1}^d \log |z \wedge w_i| \ ,
\end{equation*}
which itself is equivalent to
\begin{equation}
\label{eq:PullbackId3}
|F(z) \wedge w| = |\Res(F)|^{1/d} \prod_{i=1}^d |z \wedge w_i| \ .
\end{equation}

We verify (\ref{eq:PullbackId3}) by an explicit calculation
(compare with \cite{DeMarco}, Lemma 6.5).

Write $F = (F_1,F_2) = (\prod_{i=1}^d z \wedge a_i, \prod_{j=1}^d z \wedge b_j)$.
Since $F(b_j) = (\prod_{i=1}^d a_i \wedge b_j, 0)$ and $w = (1,0)$, 
we may take
\[
w_j = \frac{b_j}{(\prod_i a_i \wedge b_j)^{1/d}}
\]
for $j=1,\ldots,d$, where for each $j$ we fix some choice of a $d$th root of 
$\prod_i a_i \wedge b_j$.  Note that $|w_j|$ is independent of which $d$th root we
pick.  Thus
\[
\prod_{j=1}^d |z \wedge w_j| = \frac{\prod_j |z \wedge b_j|}{\prod_{i,j} |a_i \wedge b_j|^{1/d}}
= \frac{|F(z) \wedge w|}{|\Res(F)|^{1/d}} \ ,
\]
which gives (\ref{eq:PullbackId3}).
\end{proof}

Suppose now that $L$ is a complete and algebraically closed {\em non-archi\-me\-dean} field,
and let $\PP^1_{\Berk}$ denote the Berkovich projective line over $L$.
As a function on $\PP^1_{\Berk} \times \PP^1_{\Berk}$,
$g_\varphi(z,w)$ is lower semicontinuous, continuous off the diagonal, 
and symmetric in $z$ and $w$.  Also, for each fixed $z \in \PP^1_{\Berk}$,
the function $w \mapsto g_{\varphi}(z,w)$ is continuous on $\PP^1_{\Berk} \backslash \{ z \}$
(and is continuous and real-valued on all of $\PP^1_{\Berk}$ if $z \not\in \PP^1(L)$).
Moreover, $\varphi : \PP^1_{\Berk} \to \PP^1_{\Berk}$ is continuous, and both $\PP^1(L)$ 
and $\PP^1_{\Berk} \backslash \PP^1(L)$ are dense in $\PP^1_{\Berk}$.
(See \cite{BRBook}, \S7, for proofs of all of these assertions.)

In addition, it is shown in \cite{BRBook}, \S7, that for each $y \in \PP^1_{\Berk}$, there
is a natural way to define the pullback $\varphi^*(y) = \sum m_i x_i$, with $\varphi(x_i)=y$ and
$1\leq m_i \leq d$ the {\em multiplicity} of $x_i$, so that $\sum m_i = d$ and
for any continuous function $f : \PP^1_{\Berk} \to \RR$, the function $\varphi_*(f)$ defined by
\[
\varphi_*(f)(y) = \sum_{\varphi(x)=y} m_x f(x)
\]
is again continuous.
For $y \in \PP^1(L)$, $m_i$ is just the usual algebraic multiplicity of $x_i$ as
a solution to $\varphi(\cdot) = y$.

Using these facts, we can deduce that Theorem~\ref{PullbackTheorem} remains valid for all
$x,y \in \PP^1_{\Berk}$:

\begin{cor}
\label{PullbackCor}
Let $L$ be a complete and algebraically closed non-archi\-me\-dean field, and let $\varphi \in L(T)$ be a 
rational map of degree $d \geq 2$.  
Then for all $x,y \in \PP^1_{\Berk}$, we have
\begin{equation}
\label{PullbackIdentity}
g_\varphi(\varphi(x),y) = g_\varphi(x,\varphi^{*}(y)) \ .
\end{equation}
\end{cor}

\begin{proof}
Fix $y \in \PP^1(L)$.  Then since $\varphi$ is continuous, both
$g_\varphi(\varphi(x),y)$ and $g_\varphi(x,\varphi^{*}(y))$ are continuous real-valued 
functions of $x$ for all $x \in \PP^1_{\Berk} \backslash \PP^1(L)$.
Since these functions agree on the dense subset $\PP^1(L)$ of $\PP^1_{\Berk}$
by Theorem~\ref{PullbackTheorem}, 
(\ref{PullbackIdentity}) holds 
for all $x \in \PP^1_{\Berk}, y \in \PP^1(L)$.

Now fix $x \in \PP^1_{\Berk} \backslash \PP^1(L)$.  Then
$g_\varphi(\varphi(x),y)$ and $\varphi_* g_\varphi(x,y)$ are continuous and 
real-valued for all $y \in \PP^1_{\Berk}$, and they agree for $y \in \PP^1(L)$.
Therefore (\ref{PullbackIdentity}) holds 
for all $x \in \PP^1_{\Berk} \backslash \PP^1(L), y \in \PP^1_{\Berk}$.

Finally, fix $y \in \PP^1_{\Berk} \backslash \PP^1(L)$.
Then $g_\varphi(\varphi(x),y)$ and $g_\varphi(x,\varphi^{*}(y))$ are continuous and real-valued
for all $x \in \PP^1_{\Berk}$, and they agree on the dense subset
$\PP^1_{\Berk} \backslash \PP^1(L)$, so they agree for all 
$x \in \PP^1_{\Berk}$.  
\end{proof}

Using Corollary~\ref{PullbackCor}, we can deduce that 
$g_{\varphi}$ is a {\em normalized} Arakelov-Green's function
on $\PP^1_{\Berk}$.
The proof remains valid for any complete and algebraically closed field $L$ if 
we define $\PP^1_{\Berk}$ to be $\PP^1(L)$ when $L$ is archimedean.

\begin{cor}
\label{NormalizedCor}
For all $x \in \PP^1_{\Berk}$, we have
\[
\iint_{\PP^1_{\Berk} \times \PP^1_{\Berk}} 
g_{\varphi}(x,y) \mu_{\varphi}(x) \mu_{\varphi}(y) = 0 \ .
\]
\end{cor}

\begin{proof}
By (\ref{eq:GreenPotentialConstant}),
there is a constant $C$ such that
$\int g_{\varphi}(x,y) \, \mu_{\varphi}(y) = C$
for all $x \in \PP^1_{\Berk}$,
and it suffices to show that $C=0$.

Fix $x \in \PP^1_{\Berk}$.
Since $\varphi_* \mu_{\varphi} = \mu_{\varphi}$
(see \cite{BRBook}, Theorem 7.14),
it follows from Corollary~\ref{PullbackCor} that
\[
\begin{aligned}
C & = \int g_{\varphi}(x,y) \, \mu_{\varphi}(y) 
= \int g_{\varphi}(x,y) \, \varphi_* \mu_{\varphi}(y) \\
& = \int g_{\varphi}(x,\varphi(y)) \, \mu_{\varphi}(y) \\
& = \int g_{\varphi}(\varphi^*(x),y) \, \mu_{\varphi}(y) \\
& = d C \ ,  \\
\end{aligned}
\]
and thus $C=0$ as desired.
\end{proof}

Finally, we deduce a formula for the homogeneous transfinite diameter of
homogeneous filled Julia sets:

\begin{cor}
\label{TransfiniteDiameterCor}
Let $F : L^2 \to L^2$ be a non-degenerate homogeneous polynomial map,
where $L$ is a complete and algebraically closed valued field.
Let $K_F$ be the homogeneous filled Julia set of $F$.
Then
\[
d^0_\infty(K_F) = |\Res(F)|^{-\frac{1}{d(d-1)}} \ .
\]
\end{cor}

\begin{proof}
(Compare with \cite{BREQUI}, Theorem 3.16, and \cite{DeMarco}, Theorem 1.5.)
Let $R = |\Res(F)|^{-\frac{1}{d(d-1)}}$, and define
\[
D_n = \inf_{P_1,\ldots,P_n \in \PP^1(L)} \frac{1}{n(n-1)} \sum_{i \neq j} g_{\varphi}(P_i,P_j) \ .
\]

By Theorem 3.49 of \cite{BREQUI}, 
\[
\lim_{n \to \infty} D_n = \iint g_{\varphi}(x,y) \mu_{\varphi}(x) \mu_{\varphi}(y) \ ,
\]
which is zero by Corollary~\ref{NormalizedCor}.

By (\ref{eq:logscale2}), for every $P \in \PP^1(L)$ and each $\varepsilon > 0$,
there exists a lifting $z$ of $P$ to $L^2 \backslash \{ 0 \}$ 
such that $-\varepsilon \leq \Hhat_F(z) \leq 0$ (and in particular, $z \in K_F$).

Using this observation, it follows from the definitions that for every $n\geq 2$ and every 
$\varepsilon > 0$,
\[
-\log d_n^0(K_F) + \log R - 2 \varepsilon \leq D_n \leq -\log d_n^0(K_F) + \log R \ .
\]
Letting $\varepsilon \to 0$ and $n \to \infty$ gives the desired result.
\end{proof}

\section{The Mordell-Weil theorem for elliptic curves over function fields}
\label{Appendix2}

In this section, we give a new proof of 
the Mordell-Weil theorem for elliptic curves over function fields.
In addition to its simplicity and quantitative nature, our proof illustrates many of the salient features in the
proof of Theorem~\ref{maintheorem} without the complications introduced by 
the theory of Berkovich spaces.

Let $K$ be the function field of a nonsingular, integral projective algebraic curve $C$ over an algebraically closed field $k_0$,
so that $k_0$ is the field of constants of $K$.  

Let $E/K$ be an elliptic curve.
We will say that $E$ is {\em isotrivial over K}
if there is an elliptic curve $E_0/k_0$ and a $K$-isomorphism from $E$ to $E_0$,
and {\em isotrivial} if there exists a 
finite extension $K' / K$ and an elliptic curve $E'$ which is isotrivial over $K'$ 
such that $E$ is isomorphic to $E'$ over $K'$.


The following result should be compared with Theorem~\ref{IsotrivialityTheorem}.

\begin{lemma}
\label{ECLemma1}
$E$ is isotrivial if and only if $E$ has potentially good reduction
at all $v \in M_K$.
\end{lemma}

\begin{proof}
Let $j_E \in K$ denote the $j$-invariant of $E$.
By \cite{SilvermanAEC}, Proposition VIII.5.5, $E$ has potentially good reduction at $v$ if and only
if $|j_E|_v \leq 1$.  On the other hand, by  
\cite{SilvermanAEC}, Proposition III.1.4,
we see that $E$ is isotrivial if and only if 
$j_E \in k_0$.  The result now follows from the definition of $k_0$.
\end{proof}


Let $\hhat_E = \hhat_{E,K} : E(K) \to \RR$ be the N{\'e}ron-Tate canonical height on $E$.
(If we fix a Weierstrass equation $y^2 = f(x)$ for $E$ over $K$ and let
$\varphi$ be the degree four rational map expressing $x(2P)$ in terms of $x(P)$, 
then $\hhat$ can be defined by the formula $\hhat_E(P) = \frac{1}{2} \hhat_{\varphi}(x(P))$.)

According to N{\'e}ron's theory (see \cite{SilvermanAEC2}, Chapter VI, and 
\cite{SerreLMW}, \S6.5), the global canonical height
$\hhat_E$ can be decomposed into a sum of canonical local heights: for
$P \in E(K), P \neq 0$, we have
\begin{equation}
\label{eq:SumOfLocalHeights}
\hhat_E(P) = \sum_{v \in M_K} \lambda_v(P) \ ,
\end{equation}
where the functions $\lambda_v : E(\CC_v) \backslash \{ 0 \} \to \RR$ satisfy

\begin{itemize}
\item[(1)] If $E$ has good reduction at $v$, then $\lambda_v(P) \geq 0$ for all $P \in E(\CC_v) \backslash \{ 0 \}$.
\item[(2)] If $E$ has split multiplicative reduction at $v$, then 
$\lambda_v(P) = i(P) + j(P)$ for all $P \in E(\CC_v)\backslash \{ 0 \}$, where $i(P) \geq 0$ and 
$j(P) = \frac{1}{2}v(q)\B_2(r(P))$.  Here $v(q) = -\log |j_E|_v > 0$ (for an appropriate choice of logarithm), 
$r : E(\CC_v) \to \RR / \ZZ$ is the ``retraction homomorphism'' described in \cite{BakerPetsche}, \S3.1, and 
$\B_2(t) : \RR / \ZZ \to \RR$ is the periodic second Bernoulli polynomial defined for $0 \leq t \leq 1$ by
$\B_2(t) = t^2 - t + 1/6$.
\end{itemize}

By the parallelogram law, 
\begin{equation}
\label{eq:MW1}
\hhat_E(P-Q) = 2\hhat_E(P) + 2\hhat_E(Q) - \hhat_E(P+Q) \leq 2\hhat_E(P) + 2\hhat_E(Q) \ .
\end{equation}

Also, since $\frac{1}{2} \B_2(0) = \frac{1}{12} > 0$, it follows by continuity that the circle $\RR / \ZZ$ can be
decomposed as a union of
finitely many segments $U_1,\ldots,U_s$ such that $\B_2(x-y) > 0$ whenever $x,y \in U_i$ for some $i$.  For example, we can take
$U_i = [\frac{i-1}{6},\frac{i}{6}]$ for $i = 1,\ldots,6$, in which case 
\begin{equation}
\label{eq:MW2}
\frac{1}{2} \B_2(x-y) \geq \frac{1}{72}
\end{equation}
for $x,y \in U_i, \, 1\leq i \leq 6$.

Finally, a simple Fourier averaging argument (\cite{SilvermanAEC2}, Exercise 6.11) shows that
if $x_1,\ldots,x_N \in \RR / \ZZ$, then
\begin{equation}
\label{eq:MW3}
\sum_{\substack{1\leq i,j \leq N \\ i \neq j}} \frac{1}{2} \B_2(x_i - x_j) \geq -\frac{N}{12} \ .
\end{equation}

We now prove:

\begin{theorem}
\label{EllipticGapTheorem}
If $E/K$ is not isotrivial, then
there is a constant $\varepsilon > 0$ depending on $E$ and $K$ such that 
the set $\{ P \in E(K) \; : \; \hhat_E(P) \leq \varepsilon \}$ is finite.
\end{theorem}

\begin{proof}
Replacing $K$ by a finite extension, we may assume without loss of generality that
$E$ has semistable reduction over $K$, and thus that $E$ has either good or split multiplicative
reduction at each $v \in M_K$.  By Lemma~\ref{ECLemma1}, the set $S$ of places of bad reduction for $E$ is non-empty;
let $s$ be its cardinality.  For each $v \in S$, let $\delta_v = -\frac{1}{12} \log |j_E|_{v} > 0$, 
let $\delta = \min_{v \in S} \{ \delta_v \}$, and let $v_0 \in S$ be a place for which $\delta_{v_0} = \delta$.
Finally, let $S' = S \backslash \{ v_0 \}$.

Fix $N\geq 2$, let $M = 6(N-1)+1$, and suppose that $P_1,\ldots,P_M \in E(K)$ are distinct points each having 
$\hhat_E(P_i) < \delta / 48$.  

Iterating (\ref{eq:MW1}) and letting
\[
g_v(P_1,\ldots,P_N) = \frac{1}{N(N-1)} \sum_{\substack{1\leq i,j \leq N \\ i \neq j}} \lambda_{v}(P_i - P_j) \ ,
\]
we obtain
\begin{equation}
\label{eq:MW4}
\sum_{v \in M_K} g_v(P_1,\ldots,P_N)
\leq \frac{4}{N} \sum_{i} \hhat_E(P_i) 
< \frac{\delta}{12} \ . 
\end{equation}

On the other hand, by (\ref{eq:MW2}) and the pigeonhole principle, after relabelling the $P_i$'s we may assume that
$\lambda_{v_0}(P_i - P_j) \geq \delta/6$ for all $1 \leq i,j \leq N$.
Thus
\[
g_{v_0}(P_1,\ldots,P_N) \geq \delta/6 \ .
\]

In addition, for $v \in S'$, (\ref{eq:MW3}) gives
\[
g_v(P_1,\ldots,P_N) \geq -\frac{\delta_v}{N-1} \ .
\]

If $N \geq 12(s-1) + 1$, it follows that
\[
\begin{aligned}
\sum_{v \in M_K} g_v(P_1,\ldots,P_N) &\geq 
g_{v_0}(P_1,\ldots,P_N) + 
\sum_{v \in S'} g_v(P_1,\ldots,P_N) \\
&\geq \frac{\delta}{6} - \frac{\delta}{12} \\ 
& = \frac{\delta}{12} \ , \\
\end{aligned}
\]
contradicting (\ref{eq:MW4}).

Thus $N \leq \max \{ 2, 12(s-1) \}$, and therefore
$M \leq \max \{ 7, 72(s-1) - 5 \}$.  In other words, there are at most
$\max \{ 7, 72(s-1) - 5 \}$ distinct points in $E(K)$ with height less than $\delta / 48$.
\end{proof}

We recall from 
\cite{LangDG}, Proposition 6.2.3, (see also \cite{SilvermanAEC2}, Theorem III.2.1, when ${\rm char}(K)=0$)
the following ``weak Mordell-Weil theorem'' over $K$; its proof is
more or less the same as in the case where $K$ is a number field.

\begin{theorem}
\label{WeakMWThm}
If $E/K$ is not isotrivial and $E[m](\Kbar) \subseteq E(K)$, then for every integer $m \geq 1$ not divisible
by the characteristic of $K$, the group $E(K)/mE(K)$ is finite.
\end{theorem}

According to Lemma~\ref{NorthcottFailsLemma}, for $M$ sufficiently large there are
infinitely many points $x \in \PP^1(K)$ such that $\hhat_{\varphi}(x) \leq M$.
However, most of these points do not lift to $K$-rational points of $E$, 
as the following result shows:

\begin{theorem}
\label{EllipticHeightThm}
If $E/K$ is not isotrivial, then for every $M> 0$, the set
\[
\{ P \in E(K) \; : \; \hhat_E(P) \leq M \}
\]
is finite. 
\end{theorem}

\begin{proof}
(Compare with the proof of \cite{Ghioca}, Lemma 2.2.)

Let $m \geq 2$ be any integer not divisible by ${\rm char}(K)$.
Without loss of generality, we may replace $K$ by a finite extension and assume that
$E[m](\Kbar) \subseteq E(K)$.
Now suppose the result in question is false, and let
\[
C = \inf \{ M \; : \; \exists \textrm{ infinitely many } P \in E(K) \textrm{ with } \hhat_E(P) \leq M \} \ .
\]
By Theorem~\ref{EllipticGapTheorem}, $C>0$.
By the definition of $C$, there exists an infinite sequence $P_n \in E(K)$ such that $\hhat_E(P_n) < 3C/2$
for all $n$. 
Since $E(K)/mE(K)$ is finite by Theorem~\ref{WeakMWThm}, there is a coset of $mE(K)$ in $E(K)$ containing infinitely
many of the points $P_n$.  But if $i \neq j$ and $P_i$ and $P_j$ are in the same coset, then
$P_i - P_j = mQ$ for some $Q \in E(K)$, and it follows from the parallelogram law that
\[
\hhat_E(Q) \leq \frac{\hhat_E(P_i) + \hhat_E(P_j)}{m^2} < \frac{3C}{4} < C \ .
\]
It follows that there are infinitely many $Q \in E(K)$ such that $\hhat_E(Q) < \frac{3C}{4}$,
contradicting our choice of $C$.
\end{proof}

Using Theorems~\ref{WeakMWThm} and \ref{EllipticHeightThm},
the usual descent argument (see \cite{SilvermanAEC}, Proposition VIII.3.1) now establishes the following version of the Mordell-Weil theorem:

\begin{cor}[Mordell-Weil Theorem]
\label{MWCor}
If $E/K$ is not isotrivial, then $E(K)$ is finitely generated.
\end{cor}

\begin{remark}
$\left.\right.$

(i) Theorem~\ref{EllipticHeightThm} and Corollary~\ref{MWCor} are proved in \cite{LangDG}, Chapter 5,
in any characteristic under the weaker hypothesis that $E$ is not isotrivial {\em over K}.

(ii) Unlike the proofs in \cite{LangDG} or \cite{SilvermanAEC2},
our proof of Theorem~\ref{EllipticHeightThm} remains valid over an arbitrary function field $K$ 
(in the general sense discussed in the introduction).

\end{remark}


Finally, we sketch how one can deduce Theorem~\ref{EllipticGapTheorem} directly from
Theorem~\ref{maintheorem}.  For simplicity, we assume that $K$ has characteristic 
different from 2,3.  

Fix a minimal Weierstrass equation $y^2 = x^3 + Ax + B$ for $E$ over $K_v$ with discriminant $\Delta$, and
recall that the N{\'e}ron-Tate canonical height 
$\hhat_E = \hhat_{E,K} : E(K) \to \RR$ on $E$
can be defined for $P = (x(P),y(P)) \in E(K)$ as
\[
\hhat_E(P) = \frac{1}{2} \hhat_{\varphi}(x(P))
\]
where $\varphi$ is a degree $4$ rational function defined over $K$.

If $x(P) = (z_0,z_1) \in \PP^1(K)$, we have 
$\varphi(x(P)) = x([2]P) = (F_1(z_0,z_1):F_2(z_0,z_1))$ with $F_1(X,Z),F_2(X,Z) \in \O_v[X,Z]$ homogeneous
of degree 4.  
More concretely, according to \cite{SilvermanAEC}, \S~III.2.3d, we have
\[
\begin{array}{ll}
F_1(X,Z) = X^4 - 2AX^2Z^2 - 8BX + A^2 \ , \\
F_2(X,Z) = 4X^3Z + 4AXZ^3 + 4BZ^4 \ . \\
\end{array}
\]
An explicit computation shows that $\Res(F_1(X,Z),F_2(X,Z)) = 2^8 \Delta^2$.
In particular, since $K_v$ has residue characteristic different from $2$, we see that 
\begin{equation}
\label{eq:ResultantDiscriminant}
|\Res(F_1,F_2)|_v = |\Delta|^2_v \ .
\end{equation}

It is not hard to show that for all $P,Q \in E(\CC_v)$, $\lambda_v$ and $g_{\varphi,v}$ are
related by the formula
\begin{equation}
\label{eq:ComparisonFormula}
\lambda_v(P+Q) + \lambda_v(P-Q) = g_{\varphi,v}(x(P),x(Q)) \ .
\end{equation}
[{\bf Sketch of proof:} One checks that for fixed $w$, $g_{\varphi,v}(z,w)$ is a Call-Silverman canonical local height function on 
$\PP^1$ relative to the divisor $(w)$; this implies by functoriality of canonical local heights that both
sides of (\ref{eq:ComparisonFormula}) are equal up to a constant.  Finally, one shows that the constant is zero using
(\ref{eq:ResultantDiscriminant}).]

In particular, 
\begin{equation}
\label{eq:ComparisonFormula2}
\lambda_v(P) = \frac{1}{2}g_{\varphi,v}(x(P),\infty) \ .
\end{equation}

\begin{lemma}
\label{ECLemma2}
$E$ has potentially good reduction at $v \in M_K$ if and only if $\varphi$ does.
\end{lemma}

\begin{proof}
Replacing $K$ by a finite extension, we may assume without loss of generality that
$E$ has semistable reduction over $K$.
The explicit formulas for $\lambda_v$ then show that $\lambda_v(P) \geq 0$ for all
$P \in E(\CC_v)$ if and only if $E$ has good reduction at $v$.
On the other hand, Theorem~\ref{ReductionTheorem} shows that $g_v(z,w) \geq 0$
for all $z,w \in \PP^1(\CC_v)$ if and only if $\varphi$ has potential good reduction
at $v$.  The result now follows easily from (\ref{eq:ComparisonFormula}) 
and (\ref{eq:ComparisonFormula2}), together with the fact 
that the map $x : E(\CC_v) \to \PP^1(\CC_v)$ is surjective.
\end{proof}

Theorem~\ref{EllipticGapTheorem} is now an immediate consequence of Theorem~\ref{maintheorem}
and Lemmas~\ref{ECLemma1} and \ref{ECLemma2}.



\bibliographystyle{plain}
\bibliography{FFHeights}

\def\cprime{$'$} \def\cprime{$'$} \def\cprime{$'$}
\begin{thebibliography}{10}

\bibitem{BakerDGF}
M.~Baker.
\newblock A lower bound for average values of dynamical {Green's} functions.
\newblock {\em Math. Res. Lett.}
\newblock to appear. Available at {\tt arxiv:math.NT/0507484}, 15 pages.

\bibitem{BakerHsia}
M.~Baker and L.~C. Hsia.
\newblock Canonical heights, transfinite diameters, and polynomial dynamics.
\newblock {\em J. Reine Angew. Math.}, 585:61--92, 2005.

\bibitem{BakerPetsche}
M.~Baker and C.~Petsche.
\newblock Global discrepancy and small points on elliptic curves.
\newblock {\em \em Internat. Math. Res. Notices}, 2006.
\newblock to appear. Available at {\tt arxiv:math.NT/0507228}, 33 pages.

\bibitem{BRBook}
M.~Baker and R.~Rumely.
\newblock Analysis and dynamics on the {Berkovich} projective line.
\newblock preprint. Available at {\tt arxiv:math.NT/0407433}, 150 pages, 2004.

\bibitem{BREQUI}
M.~Baker and R.~Rumely.
\newblock Equidistribution of small points, rational dynamics, and potential
  theory.
\newblock {\em Ann. Inst. Fourier (Grenoble)}, 2005.
\newblock to appear. Available at {\tt arxiv:math.NT/0407426}, 50 pages.

\bibitem{BRBook2}
M.~Baker and R.~Rumely.
\newblock Potential theory on the {Berkovich} projective line.
\newblock in preparation, 180 pages, 2006.

\bibitem{BenedettoRDJ}
R.~L. Benedetto.
\newblock Reduction, dynamics, and {J}ulia sets of rational functions.
\newblock {\em J. Number Theory}, 86(2):175--195, 2001.

\bibitem{BenedettoFFHeights}
R.~L. Benedetto.
\newblock Heights and preperiodic points of polynomials over function fields,
  2005.
\newblock preprint. Available at {\tt arxiv:math.NT/0510444}.

\bibitem{BombieriGubler}
E.~Bombieri and W.~Gubler.
\newblock {\em Heights in Diophantine Geometry}, volume~3 of {\em New
  Mathematical Monographs}.
\newblock Cambridge Univ. Press, Cambridge, 2006.

\bibitem{BourbakiGT}
N.~Bourbaki.
\newblock {\em General topology. {C}hapters 1--4}.
\newblock Elements of Mathematics (Berlin). Springer-Verlag, Berlin, 1998.
\newblock Translated from the French, Reprint of the 1989 English translation.

\bibitem{DeMarco}
L.~DeMarco.
\newblock Dynamics of rational maps: {L}yapunov exponents, bifurcations, and
  capacity.
\newblock {\em Math. Ann.}, 326(1):43--73, 2003.

\bibitem{FJBook}
C.~Favre and M.~Jonsson.
\newblock {\em The valuative tree}, volume 1853 of {\em Lecture Notes in
  Mathematics}.
\newblock Springer-Verlag, Berlin, 2004.

\bibitem{Ghioca}
D.~Ghioca.
\newblock {Elliptic curves over the perfect closure of a function field}.
\newblock preprint. Available at {\tt arXiv:math.NT/0505511}, 5 pages, 2005.

\bibitem{LangDG}
S.~Lang.
\newblock {\em Fundamentals of {D}iophantine geometry}.
\newblock Springer-Verlag, New York, 1983.

\bibitem{LangAlgebra}
S.~Lang.
\newblock {\em Algebra}, volume 211 of {\em Graduate Texts in Mathematics}.
\newblock Springer-Verlag, New York, third edition, 2002.

\bibitem{Moriwaki}
Atsushi Moriwaki.
\newblock Arithmetic height functions over finitely generated fields.
\newblock {\em Invent. Math.}, 140(1):101--142, 2000.

\bibitem{RLPHS}
J.~Rivera-Letelier.
\newblock Espace hyperbolique {$p$}-adique et dynamique des fonctions
  rationnelles.
\newblock {\em Compositio Math.}, 138(2):199--231, 2003.

\bibitem{SerreLMW}
J.-P. Serre.
\newblock {\em Lectures on the {M}ordell-{W}eil theorem}.
\newblock Aspects of Mathematics. Friedr. Vieweg \& Sohn, Braunschweig, third
  edition, 1997.
\newblock Translated from the French and edited by Martin Brown from notes by
  Michel Waldschmidt, With a foreword by Brown and Serre.

\bibitem{SilvermanAEC}
J.~H. Silverman.
\newblock {\em The arithmetic of elliptic curves}, volume 106 of {\em Graduate
  Texts in Mathematics}.
\newblock Springer-Verlag, New York, 1986.

\bibitem{SilvermanAEC2}
J.~H. Silverman.
\newblock {\em Advanced topics in the arithmetic of elliptic curves}, volume
  151 of {\em Graduate Texts in Mathematics}.
\newblock Springer-Verlag, New York, 1994.

\bibitem{Thuillier}
A.~Thuillier.
\newblock {\em Th{\'e}orie du potentiel sur les courbes en g{\'e}om{\'e}trie
  analytique non archim{\'e}dienne. {A}pplications {\`a} la th{\'e}orie
  d'{A}rakelov}.
\newblock PhD thesis, University of Rennes, 2005.
\newblock preprint. Available at {\tt
  http://tel.ccsd.cnrs.fr/documents/archives0/00/01/09/90/index.html}.

\end{thebibliography}


\end{document}